\documentclass[11pt]{article}

\usepackage{latexsym,amsmath,amsfonts,amssymb,amstext,theorem,euscript,array,enumerate,mathrsfs}

\newtheorem{Theorem}{Theorem}[section]
\newtheorem{Definition}{Definition}[section]
\newtheorem{Proposition}{Proposition}[section]

\newtheorem{Lemma}{Lemma}[section]

\newtheorem{Remark}{Remark}[section]
\newtheorem{Example}{Example}[section]
\numberwithin{equation}{section}

\def\Dim{\noindent\hbox{{\bf Proof.}$\;\; $}}
\def\finedim{{\hfill\hbox{\enspace${ \square}$}} \smallskip}

\def \X{\mathfrak{X}}
\def \I{\mathfrak{I}}
\def \N{\mathbb{N}}
\def \R{\mathbb{R}}

\def \E{\mathbb{E}}

\def \P{\mathbb{P}}

\def \Ac{{\cal A}}
\def \Acn{{\cal A}_n}

\def \Fc{{\cal F}}

\def \Nc{{\cal N}}
\def \Pc{{\cal P}}

\def \Rc{{\cal R}}

\def \Uc{{\cal U}}

\def \ep{\hbox{ }\hfill$\Box$}

\allowdisplaybreaks

\begin{document}

\title{Ergodic Control of Infinite Dimensional SDEs\\ with Degenerate Noise}

\author{ Andrea COSSO\footnote{Department of Mathematics, University of Bologna, Piazza di Porta S. Donato 5, 40126 Bologna,
Italy, {\tt andrea.cosso@unibo.it}.}
\qquad\quad
Giuseppina GUATTERI\footnote{Department of Mathematics, Politecnico di Milano, Piazza Leonardo da Vinci 32, 20133 Milano, Italy, {\tt giuseppina.guatteri@polimi.it}.}
\qquad\quad
Gianmario TESSITORE\footnote{Department of Mathematics and Applications, University of Milano-Bicocca, Via R. Cozzi 53 - Building U5, 20125 Milano, Italy, {\tt gianmario.tessitore@unimib.it}.}} 

\maketitle

\begin{abstract}
The present paper is devoted to the study of the asymptotic behavior of the value functions of both finite and infinite horizon stochastic control problems and to the investigation of their relation with suitable stochastic ergodic control problems. Our methodology is based only on probabilistic techniques, as for instance the so-called randomization of the control method, thus avoiding completely analytical tools from the theory of viscosity solutions. We are then able to treat with the case where the state process takes values in a general (possibly infinite dimensional) real separable Hilbert space and the diffusion coefficient is allowed to be degenerate.
\end{abstract}

\vspace{5mm}


\vspace{5mm}


\date{}

\section{Introduction}

In the present paper we study the asymptotic behavior  of the value functions both for  finite horizon stochastic control problems (as the horizon diverges) and for  discounted  infinite horizon control problems (as the discount vanishes) and investigate their relation with suitable stochastic ergodic control problems. We will refer to such limits as \textit{ergodic limits}. The main novelty of this work is that we deal with ergodic  limits for  control problems in which the state process is allowed to take values in a general (possibly infinite dimensional) real separable Hilbert space and the   diffusion coefficient is allowed to be degenerate.

 On the one side, the infinite dimensional framework imposes the use of purely probabilistic techniques essentially based on backward stochastic differential equations (BSDEs for short, see for instance the introduction of \cite{FuTessAOP2002}), on the other, the degeneracy of the noise prevents the use of standard BSDEs techniques as they are, for instance, implemented, for similar problems, in \cite{FuhHuTess}.  Indeed, see again \cite{FuTessAOP2002}, the identification between solutions  of BSDEs and value functions of stochastic optimal control problems can be easily obtained, by a Girsanov argument, as far as the  so called  \textit{structure condition},  imposing large enough image of diffusion operator, holds. Here we wish to avoid such a requirement.

To our knowledge, the only paper that deals, by means of BSDEs,  with ergodic limits in the  degenerate case  is \cite{CoFuhPham} where authors use the same tool of randomized control problems and related \textit{constrained} BSDEs that we will eventually employ here.
Notice however that in \cite{CoFuhPham} the state process lives in a finite-dimensional Euclidean space and probabilistic methods are combined with PDE techniques, relying on powerful tools from the theory of viscosity solutions. Here, as already mentioned, we have to completely avoid these arguments. As a matter of fact viscosity solutions require, in the infinite dimensional case, additional artificial assumptions  that we can not impose here (see for instance the theory of $B$-continuous viscosity solutions for second order PDEs, \cite{FGS}, \cite{Swiech}). On the other side to separate difficulties we consider, as in \cite{FuhHuTess} but differently from \cite{CoFuhPham}, only additive and uncontrolled noise. This in particular considerably simplifies the proof of estimate \eqref{stimaZpenunif}.

Let us now give a more precise idea of the results obtained in the paper. Consider the following infinite and finite horizon stochastic control problems:
\[
v^\beta(x) \ = \ \inf_{u\in\Uc}\E\bigg[\int_0^\infty e^{-\beta t}\,\ell(X_s^{x,u},u_s)\,ds\bigg]
\]
and
\[
v^T(x) \ = \ \inf_{u\in\Uc}\E\bigg[\int_0^T \ell(X_{s}^{x,u},u_{s})\,ds+\phi(X_{T}^{x,u})\bigg],
\]
where the discount coefficient $\beta$ can be any positive real number, as well as the time horizon $T$, while the controlled state process $X^{x,u}$ takes values in some real separable Hilbert space $H$ and is a (mild) solution to the time-homogenous stochastic differential equation
\[
dX_t \ = \ AX_tdt + F(X_t,u_t)dt + Gd W_t^1, \qquad X_0 \ = \ x.
\]
Here $W^1$ is a cylindrical Wiener process and $A$ is a possibly unbounded linear operator on $H$. We assume that  both $A$ and $F$ are dissipative. The control process $u$ is progressively measurable and takes values in some real separable Hilbert space $U$ (actually, $U$ can be taken more general, see Remark \ref{R:U_Hilbert}). Notice that the diffusion coefficient $G$ is only assumed to be a bounded linear operator, so that it can be degenerate. Our aim is to study the asymptotic behavior of
\begin{equation}\label{Asympt}
v^\beta(x)-v^\beta(0), \qquad\qquad \beta\,v^\beta(0), \qquad\qquad \frac{v^T(x)}{T},
\end{equation}
as $\beta\rightarrow0$ and $T\rightarrow+\infty$. In order to do it, we find non-linear Feynman-Kac representations for both $v^\beta$ and $v^T$ in terms of suitable backward stochastic differential equations (which can be seen as the probabilistic counterparts of the Hamilton-Jacobi-Bellman equations). Since $G$ can be degenerate we adopt the recently introduced so-called randomization method, see e.g. \cite{KharroubiPham}, \cite{FuhrmanPham}, \cite{BCFP}, which was also implemented in \cite{CoFuhPham}. Here we use it in a rather different way, as we will explain below. The idea of the randomization (of the control) method is to introduce a new control problem (called the randomized infinite/finite horizon stochastic control problem), where we replace the family of control processes $u$ by a particular class of processes (depending on a control parameter $\alpha$), here denoted by $\I^{a,\alpha}$, which is, roughly speaking, dense in $\Uc$. More precisely, focusing for simplicity only on the infinite horizon case, we define the value function of the randomized infinite horizon stochastic control problem as follows:
\[
v^{\beta,\Rc}(x,a) \ = \ \inf_{\alpha\in\Ac}\E\bigg[\int_0^\infty e^{-\beta t}\,\ell(\X_t^{x,a,\alpha},\I_t^{a,\alpha})\,dt\bigg],
\]
where $\Ac$ is the set of progressively measurable and uniformly bounded processes taking values in $U$, while the state process is the pair $(X^{x,a,\alpha},\I^{a,\alpha})$ satisfying
\[
\begin{cases}
d\X_t\ &= \ A\X_t dt + F(\X_t, \I_t)dt +\,GdW_t^1, \hspace{8.7mm} \X_0 \ = \ x, \\
d \I_t\ &= \ R\alpha_t dt + RdW_t^2,  \hspace{33.8mm} \I_0 \ = \ a, \\
\end{cases}
\]
with $R\colon U\rightarrow U$ being a trace class injective linear operator with dense image, while $W^2$ is a cylindrical Wiener process independent of $W^1$. We prove by a density argument (Proposition \ref{P:v=V^R}) that $v^\beta(x)=v^{\beta,\Rc}(x,a)$, for every $(x,a)\in H\times U$, so, in particular, $v^{\beta,\Rc}$ does not depend on its second argument $a$. Notice that we randomize the control by means of an independent cylindrical Wiener process $W^2$, instead of using an independent Poisson random measure on $\R_+\times U$ (as it is usually the case in the literature on the randomization method). Taking a Poisson random measure has the advantage that $U$ can be any Borel space, while here we have to impose some restrictions on $U$ (see Remark \ref{R:U_Hilbert}). However, randomizing the control by means of a cylindrical Wiener process is simpler  and more natural in an infinite dimensional setting, as many important fundamental results on SDEs and BSDEs can only be found for the case where the driving noise is of Wiener type. Moreover, with this choice the results presented here can receive
more attention in the infinite dimensional literature. Furthermore, the Wiener type randomization has not been enough investigated in the literature, since it was implemented only in \cite{ChoukrounCosso}, where however the proof of the fundamental equality $v^{\beta}(x)=v^{\beta,\Rc}(x,a)$ was based on PDE techniques (in particular, viscosity solutions' arguments) adapted from \cite{KharroubiPham}, instead of using purely probabilistic arguments, as it was done in \cite{FuhrmanPham} and \cite{BCFP} for the Poisson type randomization. So, in particular, this is the first time that the equality $v^{\beta}(x)=v^{\beta,\Rc}(x,a)$ is proved in purely probabilistic terms for the Wiener type randomization.

Once we know that $v^\beta(x)=v^{\beta,\Rc}(x,a)$, it is fairly standard in the framework of the randomization approach,  to derive a non-linear Feynman-Kac formula for $v^\beta$, see for instance \cite{KharroubiPham}. As a matter of fact, notice that, for each positive integer $n$, the control problem with value function
\[
v^{\beta,\Rc,n}(x,a) \ = \ \inf_{\alpha\in\Ac: |\alpha| \leq n }\E\bigg[\int_0^\infty e^{-\beta t}\,\ell(X_t^{x,a,\alpha},I_t^{a,\alpha})\,dt\bigg]
\]
is a dominated problem. Therefore, by standard BSDE techniques, $v^{\beta,\Rc,n}$ admits a non-linear Feynman-Kac representation in terms of some BSDE depending on the parameter $n$. Passing to the limit as $n$ goes to infinity, we  find, as in \cite{CoFuhPham}, a non-standard BSDE for $v^\beta(x)=v^{\beta,\Rc}(x,a)=\lim_n v^{\beta,\Rc,n}(x,a)$ (see Propositions \ref{PenEllRand} and \ref{propEllitica}). As a matter of fact such a BSDE, that we shall call `constrained', involves a reflection term and is characterized by its maximality, see again \cite{KharroubiPham}. We eventually exploit the probabilistic representation of $v^\beta$ (and similarly of $v^T$) to study the limits in \eqref{Asympt}. In particular, in Section \ref{S:Ergodic} we prove that, up to a subsequence, the limit of $v^\beta(x)-v^\beta(0)$ (resp. $\beta v^\beta(0)$) exists and is given by a function $\hat v$ (resp. constant $\lambda$). Moreover we prove that $\hat v$ and $\lambda$  are related to a suitable constrained ergodic backward stochastic differential equation, again of the non standard, constrained, type see Theorem \ref{ergodico}. This is the the most technical  result of the paper. In the previous literature, see \cite{CoFuhPham}, PDE techniques are indeed  of great help at this level. In  the present context we have to prove in a direct way that the candidate solution to the ergodic constrained BSDE enjoys the required maximality property. To do that we exploit the extra regularity of the trajectories of the state equation implied by  Assumption ({\bf A.2}) (see estimate \eqref{stima_compattezza} in Proposition \ref{Eq-stato.randomizzata}).

Concerning the long time asymptotics of $v^T(x)/T$, we show in Theorem \ref{T:LongTimeAsympt} that this quantity converges to the same constant $\lambda$. We end Section \ref{S:Ergodic} proving that, under suitable assumptions, $\lambda$ coincides with the value function of an ergodic control problem. This latter result is again proved using only probabilistic techniques, while in \cite{CoFuhPham} the proof is based on PDE arguments (see Remark \ref{R:lambda} for more details on this point).

The rest of the paper is organized as follows. In Section \ref{S:InfHor} we firstly introduce the notations used throughout the paper, then we formulate both the infinite and finite horizon stochastic optimal control problems on a generic probabilistic setting; afterwards, we formulate both control problems on a specific probabilistic, product-space, setting and we prove (Proposition \ref{UguaglianzaValueFunction1}) that, even if the probabilistic setting has changed, value functions are still the same. Section \ref{S:RandInfHor} is devoted to the formulation of the randomized control problems; we prove (Proposition \ref{P:v=V^R}) that the value functions of the randomized problems coincide with the value functions of the original control problems. In Section \ref{S:BSDE} we find non-linear Feynman-Kac representation formulae for the value functions in terms of constrained backward stochastic differential equations. In Section \ref{S:Ergodic} we introduce an ergodic BSDE and study the asymptotic behavior of the infinite horizon problem. Finally, in Section \ref{S:Finite} we introduce an ergodic control problem and study the long-time asymptotics of the finite horizon problem.

\section{Infinite/finite horizon optimal control problems}
\label{S:InfHor}

\setcounter{equation}{0} \setcounter{Assumption}{0}
\setcounter{Theorem}{0} \setcounter{Proposition}{0}
\setcounter{Corollary}{0} \setcounter{Lemma}{0}
\setcounter{Definition}{0} \setcounter{Remark}{0}

In the present section we introduce both an infinite and a finite horizon stochastic optimal control problem, first on a generic probability space and then on an enlarged probability space in product form (this latter will be then use throughout the paper). Firstly we fix some notations.

\subsection{General notation}

Let $\Xi$, $H$ and $U$ be real separable Hilbert spaces. In the sequel, we use the notations $|\cdot|_\Xi$, $|\cdot|_H$ and $|\cdot |_U$ to denote the norms on $\Xi$, $H$ and $U$ respectively; if no confusion arises, we simply write $|\cdot|$. We use similar notation for the scalar products. We denote the dual spaces of $\Xi$, $H$ and $U$ by $\Xi^*$, $H^*$,  and $U^*$ respectively. We also denote by $L(H,H)$  the space of bounded linear operators from $H$ to $H$, endowed with the operator norm. Moreover, we denote by $L_2(\Xi,H)$ the space of Hilbert-Schmidt operators from $\Xi$ to $H$. Finally, we denote by $\mathcal{B}(\Lambda)$ the Borel
$\sigma$-algebra of any topological space $\Lambda$.

Given a complete probability space $(\Omega, \mathcal{F}, \mathbb{P})$ together with a filtration $(\mathcal{F}_t)_{t\geq0}$ (satisfying the usual conditions of $\P$-completeness and right-continuity) and  an arbitrary real separable Hilbert space $V$  we define the following classes of processes for fixed $0\leq t\leq T $ and  $p\geq 1$:
\begin{itemize}
\item   $L^p_\Pc (\Omega\times [t,T];V)$ denotes the set of (equivalence classes) of $(\mathcal{F}_t)$-predictable processes $Y \in L^p (\Omega\times
[t,T];V)$ such that the following norm is finite:
\[
|Y|_p \ = \ \bigg(\E \int_t^T |Y_s|^p \, ds\bigg)^{1/p}
\]
\item $L^{p, loc}_{\mathcal{P}}(\Omega\times [0,+\infty[;V)$ denotes the set of processes defined on $\mathbb{R}^+$, whose restriction to an arbitrary time interval $[0,T]$ belongs to $L^p_\Pc (\Omega\times [0,T];V) $.  
   \item $L^p_\Pc(\Omega;C([t,T];V))$
     denotes the set of $(\mathcal{F}_t)$-predictable processes $Y$ on $[t,T]$ with continuous paths in $V$, such
    that the norm
    \[  \|Y\|_p \ = \ \big(\E \sup _{s \in [t,T]} |Y_s|^p\big)^{1/p}\]
    is finite. The elements of $L^p_{\mathcal{P}}(\Omega;C([t,T];V))$
    are identified up to indistinguishability.
    \item $L^{p, loc}_\Pc (\Omega;{C} [0,+\infty[;V)$ denotes the set of processes defined on $\mathbb{R}^+$, whose restriction to an arbitrary time interval $[0,T]$ belongs to $L^p_\Pc(\Omega;C([0,T];V))$. 
    \item $\mathcal{K}_{\mathcal{P}}^2(0,T)$ denotes the set of real-valued $(\mathcal{F}_t)$-adapted nondecreasing continuous processes $K$ on $[0,T]$ such that $\E|K_T|^2< \infty$ and $K_t=0$.
    \item $\mathcal{K}^{2,loc}_{\mathcal{P}}$ denotes the set of processes defined on $\mathbb{R}^+$, whose restriction to an arbitrary time interval $[0,T]$ belongs to $\mathcal{K}^2_{\mathcal{P}}(0,T)$.
\end{itemize}

\subsection{Formulation of the control problems}\label{sub-originalcontrol}

We formulate here both the discounted, infinite horizon, control problem and the finite horizon one whose asymptotic behavior is the main focus of the present paper. The notation chosen here may seem a bit artificial but this is done in order to keep the notation simple in the product space and randomized setting (see Sections 2 and 3)
where the technical arguments are developed.

\medskip 

We fix   a complete
probability space $(\bar\Omega^1,\bar\Fc^1,\bar\P^1)$ on which  a cylindrical Wiener process $\bar W^1 = (\bar W^1_t)_{t\geq 0}$  with
values in $\Xi$ is defined.
By $(\bar\Fc_t^1)_{t\geq0}$, or simply $(\bar\Fc_t^1)$, we denote the natural filtration
of $\bar W^1$, augmented with the family $\bar{\mathcal N}^1$ of
$\bar\P^1$-null sets of $\bar\Fc^1$. Obviously, the filtration
$(\bar\Fc_t^1)$ satisfies the usual conditions of right-continuity and $\bar\P^1$-completeness.

 In this section the notion of measurability and progressive measurability will always refer to the filtration $\bar\Fc^1$.

Let $\bar\Uc$ be the family of $(\bar\Fc_t^1)$-progressively measurable processes taking values in $U$ (see Remark \ref{R:U_Hilbert} below for the case where the space of control actions $U$ is not necessarily a Hilbert space).

\paragraph{State process.} Given $x\in H$ and $\bar u\in\bar\Uc$, we consider the controlled stochastic differential equation
\begin{equation}\label{State1}
d\bar X_t \ = \ A\bar X_tdt + F(\bar X_t,\bar u_t)dt + Gd\bar W_t^1, \qquad \bar X_0 \ = \ x.
\end{equation}
On the coefficients $A$, $F$, $G$ we impose the following assumptions.

\begin{itemize}
\item[({\bf A.1})] $A\colon\mathcal{D}(A)\subset H\to H$ is a linear, possibly unbounded operator   generating an  analytic semigroup  $\{e^{tA}\}_{t\geq 0}$. We assume that $A$ is dissipative i.e. $<Ax,x>\,\leq 0$, for all $x\in \mathcal{D}(A)$.

\item[({\bf A.2})] $G\colon\Xi\rightarrow H$ is a bounded linear operator. Moreover,
 there exist positive constants $M_A$ and $\gamma\in[0,\frac{1}{2}[$ such that
$$
\big|e^{sA}G\big|_{L_2(\Xi,H)} \ \leq \ \frac{M_A}{ s^\gamma},\qquad 
\hbox{
for all $s\in (0,1)$.}$$

\item[({\bf A.3})]  Fixed $\delta >$0  we denote by $\mathcal{D}((\delta I-A)^\rho)  $, the domain of the fractional power of the operator $ \delta I -A$, see \cite{Lunardi}. We assume that there exists a $\rho \in (0,  
\frac{1}{2}-\gamma) $ such that the domain of the fractional power $\mathcal{D}((\delta I-A)^\rho)  $ is compactly embedded in $H$.
\item[({\bf A.4})] $F\colon H\times U\rightarrow H$ is continuous and there exists  $ C_F>0$ such that
\[  
|F(x,a)| \leq C_F (1+ |x|)
\]
for all $x \in H$ and $a \in U$.

Moreover there exists $L_F>0$ such that
\[
|F(x,a) - F(x',a)|_H \ \leq \ L_F  |x - x'|_H,
\]
for all $x,x'\in H$ and $a\in U$.
\item[({\bf A.5})] $F$ is assumed to be strongly dissipative: there exists $\mu>0$ such that
\[
\langle F(x,a) -  F(x',a), x - x' \rangle_H \ \leq \ -\mu|x - x'|_H^2,
\]
for all $x, x'\in H$ and $a \in U$.

\end{itemize}
\begin{Remark} Assumption {\bf(A.4)} can be balanced  with the dissipativity  of $A$ replacing $A$ by $A+\lambda I$ and $F$ by $F-\lambda I$ ($\lambda\in \mathbb{R}$). In particular we can always think that both $A$ and $F$ are strongly dissipative.
\end{Remark}

\begin{Proposition}\label{Prop-esunstate}
Assume {\bf(A.1)}--{\bf(A.5)}. Then, for any $x\in H$ and $\bar u\in\bar\Uc$, there exists a unique (up to indistinguishability) process $\bar X^{x,\bar u}=(\bar X_t^{x,\bar u})_{t\geq0}$ that belongs to $L^{p,loc}_\Pc (\bar\Omega^1;C([0,+\infty[;H))$ for all $p \geq 1$ and is a mild solution of \eqref{State1}, that is:
\[
\bar X_t^{x,\bar u} \ = \ e^{tA} x  + \int_0^t e^{(t-s)A} F(\bar X_s^{x,\bar u},\bar u_s) \, ds +  \int_0^t e^{(t-s)A}G \, d\bar W_s^1, \qquad \text{for all }t\geq0,\,\bar\P^1\text{-a.s.}
\]
Moreover the following estimates hold:
\begin{itemize}
\item for every $T>0$ and $p\geq 1$ there exists a positive constant $\kappa_{p,T}$, independent of $x\in H$ and $\bar u\in\bar\Uc$, such that
\[
\bar{\E}^1\Big[\sup_{t \in [0,T]} |\bar{ X}^{x,\bar{u}}_t|^p\Big] \ \leq \ \kappa_{p,T} (1+ |x|^p);
\]
\item there exists a positive constant $\kappa$, independent of $x\in H$, $\bar u\in\bar\Uc$, $t\geq 0$, such that
\[
\bar{\E}^1|\bar{ X}^{x,\bar{u}}_t|\ \leq \ \kappa(1+ |x|).
\]

\end{itemize}
\end{Proposition}
\textbf{Proof.}
This is a standard result, see \cite[Proposition 3.6] {FuhHuTess} for the proof in a general Banach space context. 
Notice that the presence of the control process $\bar u$ does not causes any additional difficulty since  assumptions ({\bf A.4}) and ({\bf A.5}) hold uniformly with respect to the control variable $a\in U$.
\ep
\medskip

\noindent Finally, we fix a running cost  $\ell\colon H\times U\rightarrow\R$  and we impose the following assumption.
\begin{itemize}
\item[({\bf A.6})]
$\ell$ is continuous and bounded, moreover there exists $L_\ell>0$ such that
\[
|\ell(x,a) - \ell(x',a)| \ \leq \ L_\ell |x - x'|_H ,
\]
for all $x,x'\in H$ and $a\in U$.
\end{itemize}
\paragraph{Infinite horizon control problem.} Given a positive discount $\beta>0$, the  cost corresponding to control $\bar u\in \Uc$ and  initial condition $x$ is defined as
\[
\bar J^\beta(x,\bar u) \ := \ \bar\E^1\bigg[\int_0^\infty e^{-\beta t}\,\ell(\bar X_t^{x,\bar u},\bar u_t)\,dt\bigg],
\]
where $\bar\E^1$ denotes the expectation with respect to $\bar\P^1$. Moreover the value function is given by 
\[
\bar v^\beta(x) \ := \ \inf_{\bar u\in\bar\Uc}\bar J^\beta(x,\bar u), \qquad \text{for every }x\in H.
\]

\paragraph{Finite horizon control problem.} Fix a function $\phi: H\rightarrow \mathbb{R}$ 
satisfying:
\begin{itemize}
\item[({\bf A.7})] $\phi$ is continuous and there exists $C_\phi>0$ such that
\[
|\phi(x)| \ \leq \ C_\phi(1 + |x|_H),
\qquad\hbox{
for all $x\in H$}.\]
\end{itemize} The  cost with (finite) horizon $T>0$ and discount $\beta\geq 0$ relative to the control $\bar u$  and initial condition $x$ is defined as
\[
\bar J^{\beta,T}(x,\bar u) \ := \ \bar\E^1\bigg[\int_0^T  e^{-\beta s}\ell(\bar X_{s}^{x,\bar u},\bar u_{s})\,ds+ e^{-\beta T}\phi(\bar X_{T}^{x,\bar u})\bigg]
\]
Finally, the value function is given by 
\[
\bar v^{\beta,T}(x) \ := \ \inf_{\bar u\in\bar\Uc}\bar J^{\beta,T}(x,\bar u), \qquad \text{for every }x\in H.
\]

\begin{Remark}\label{R:U_Hilbert}
The request on the space of control actions $U$ to be an Hilbert space can be relaxed. As a matter of fact, suppose that the space of control actions is  a certain set $\tilde U$, so that, in the formulation of the stochastic optimal control problem, drift and running cost are defined on $H\times\tilde U$:
\[
\tilde F\colon H\times\tilde U\rightarrow H, \qquad\qquad \tilde\ell\colon H\times\tilde U\rightarrow\R.
\]
Suppose that $\tilde U$ has the following property: there exists a continuous surjection $\varphi\colon U\rightarrow\tilde U$, for some real separable Hilbert space $U$. This holds true, for instance, if $\tilde U$ is a compact, connected, locally connected subset of $\R^n$, for some positive integer $n$ (in this case, the existence of a continuous surjection $\varphi\colon U\rightarrow\tilde U$, with $U=\R$, follows from the Hahn-Mazurkiewicz theorem, see for instance Theorem 6.8 in \cite{Sagan}). Then, we define $F\colon H\times U\rightarrow H$ and $\ell\colon H\times U\rightarrow\R$ as
\[
F(x,a) \ := \ \tilde F(x,\varphi(a)), \qquad\qquad \ell(x,a) \ := \ \tilde\ell(x,\varphi(a)),
\]
for every $(x,a)\in H\times U$. Notice that, if $\tilde F$ (resp. $\tilde\ell$) satisfies assumptions {\bf(A.4)} and {\bf(A.5)} (resp. {\bf(A.6)}) then $F$ (resp. $\ell$) still satisfies the same assumptions. Replacing $\tilde U$, $\tilde F$, $\tilde\ell$ by $U$, $F$, $\ell$ we find a stochastic control problem of the form studied in the present work, which has the same value function of the original control problem.
\end{Remark}

\subsection{Formulation of the control problems on a product space}
\label{SubS:InfHorProductSpace}

For a technical reason imposed by the randomization method (see the next Section \ref{S:RandInfHor}) we have to reformulate our control  problems in a product probability space. The main point of this section, see Proposition \ref{UguaglianzaValueFunction1}, will be to show that this new setting does not affect the value function.

\medskip

Let $\bar W^2 = (\bar W^2_t)_{t\geq 0}$ be a cylindrical Wiener process with
values in $U$, defined on a
complete
probability space $(\bar\Omega^2,\bar\Fc^2,\bar\P^2)$. We define $(\Omega,\Fc,\P)$, $W^1$, $W^2$ as follows: $\Omega:=\bar\Omega^1\times\bar\Omega^2$, $\Fc$ the $\bar\P^1\otimes\bar\P^2$-completion of $\bar\Fc^1\otimes\bar\Fc^2$, $\P$ the extension of $\bar\P^1\otimes\bar\P^2$ to $\Fc$, $W^1(\omega^1,\omega^2):=\bar W^1(\omega^1)$, $W^2(\omega^1,\omega^2):=\bar W^2(\omega^2)$, for every $(\omega^1,\omega^2)\in\Omega$.

By $(\Fc_t)_{t\geq0}$, we denote the natural filtration
of $(W^1,W^2)$, augmented with the family $\mathcal N$ of
$\P$-null sets of $\Fc$. Clearly, $(\Fc_t)$ satisfies the usual conditions of right-continuity and $\P$-completeness. Finally, we denote by $\Uc$ the family of $(\Fc_t)$-progressively measurable processes with values in $U$. In this section measurability will always be referred to such a filtration.

\medskip As before, given $x\in H$ and $u\in\Uc$, we consider the controlled stochastic differential equation
\begin{equation}\label{State}
dX_t \ = \ AX_tdt + F(X_t,u_t)dt + Gd W_t^1, \qquad X_0 \ = \ x.
\end{equation}
Exactly as  for equation \eqref{State1}, we have the following result.
\begin{Proposition}\label{RegStatoU}
Assume {\bf(A.1)}--{\bf(A.5)}. Then, for any $x\in H$ and $u\in\Uc$, there exists a unique (up to indistinguishability) process $X^{x,u}=(X_t^{x,u})_{t\geq0}$ that  belongs to $L^{p,loc}_\Pc (\Omega;C([0,+\infty[;H))$ for all $p \geq 1$ and is a mild solution of \eqref{State}, that is:
\[
X_t^{x, u} \ = \ e^{tA} x  + \int_0^t e^{(t-s)A} F(X_s^{x, u}, u_s) \, ds +  \int_0^t e^{(t-s)A}G \, dW_s^1, \qquad \text{for all }t\geq0,\,\P\text{-a.s.}
\]
Moreover the following estimates hold:
\begin{itemize}
\item for every $T>0$ and $p\geq 1$ there exists a positive constant $\kappa_{p,T}$, independent of $x\in H$ and $u \in\Uc$, such that
\begin{equation}
\label{stimadisstato-oriz-fin-prod}
{\E}\Big[\sup_{t \in [0,T]} |{ X}^{x,{u}}_t|^p\Big] \ \leq \ \kappa_{p,T} (1+ |x|^p);
\end{equation}
\item there exists a positive constant $\kappa$, independent of $x\in H$, $u \in\Uc$, $t\geq 0$, such that
\begin{equation}
\label{stimadisstato-oriz-infin-prod}
{\E}|\bar{ X}^{x,{u}}_t|\ \leq \ \kappa(1+ |x|).
\end{equation}
\end{itemize}
\end{Proposition}

\noindent Again, for every $\beta>0$, and for any $x\in H$, $u\in\Uc$, we define the \textit{infinite horizon} cost
\[
J^\beta(x,u) \ := \ \E\bigg[\int_0^\infty e^{-\beta t}\,\ell(X_t^{x,u},u_t)\,dt\bigg]
\]
and the corresponding value function
\[
v^\beta(x) \ := \ \inf_{u\in\Uc}J^\beta(x,u), \qquad \text{for every }x\in H.
\]

\noindent On the other hand, for every $T>0$, $\beta\geq 0$, and for any $x\in H$, $u\in\Uc$, we define the \textit{finite horizon} cost:
\[
J^{\beta,T}(x,u) \ := \ \E\bigg[\int_0^T e^{-\beta s}\ell(X_{s}^{x,u},u_{s})\,ds+e^{-\beta T}\phi(X_{T}^{x,u})\bigg].
\]
and the corresponding value function 
\[
v^{\beta,T}(x) \ := \ \inf_{u\in\Uc}J^{\beta,T}(x,u).
\]

\medskip

\noindent In the next proposition we give a detailed argument showing  that, as expected, the value function of both infinite and finite horizon problems is not affected by the product space formulation. For the definition of  $\bar v^\beta$ and $\bar v^T$ see Section
\ref{sub-originalcontrol}.

\begin{Proposition}\label{UguaglianzaValueFunction1}
Suppose that Assumptions {\bf(A.1)}--{\bf(A.7)} hold. Then:
\begin{itemize}
\item[\textup{(i)}] For all $\beta>0$, $\bar v^\beta(x)=v^\beta(x)$, for every $x\in H$.
\item[\textup{(ii)}] For all $T>0$ and $\beta\geq 0$, $\bar v^{\beta,T}(x)=v^{\beta,T}(x)$, for every $x\in H$.
\end{itemize}
\end{Proposition}
\textbf{Proof.}
We prove only the first statement, since the proof of (ii) can be done proceeding along the same lines.

Fix $\beta>0$ and $x\in H$. We begin noting that the inequality $\bar v^\beta(x)\geq v^\beta(x)$ is immediate. As a matter of fact, given $\bar u\in\bar\Uc$ (thus $\bar u $ is a process defined on $[0,\infty[\times \bar \Omega^1$) let $u_t(\omega^1,\omega^2):=\bar u_t(\omega^1)$, for every $(\omega^1,\omega^2)\in\Omega$. Then, $u\in\Uc$, and $\bar J^\beta(x,\bar u)=J^\beta(x, u)\geq v^\beta(x)$. Taking the infimum over $\bar u\in\bar\Uc$, we conclude that $\bar v^\beta(x)\geq v^\beta(x)$.

We now prove the other inequality. To this end, we recall that $(\bar\Fc_t^1)$ is the natural filtration on $(\bar\Omega^1,\bar\Fc^1,\bar\P^1)$ of $\bar W^1$, augmented with the family $\bar\Nc^1$ of $\bar\P^1$-null sets. In a similar way, we define $(\bar\Fc_t^2)$. Now, let $(\tilde\Fc_t)$ be the filtration defined as $\tilde\Fc_t:=\bar\Fc_t^1\otimes\bar\Fc_t^2$, for every $t\geq0$. Observe that $(\tilde\Fc_t)$ is right-continuous (as it can be shown proceeding for instance as in the proof of Theorem 1 in \cite{HeWang}), but not necessarily $\P$-complete. We also notice that $(\Fc_t)$ is the augmentation of $(\tilde\Fc_t)$.

Now, fix $ u\in\Uc$. Since $ u$ is $(\Fc_t)$-progressively measurable, by for instance Lemma B.21 in \cite{BainCrisan} or Theorem 3.7 in \cite{ChungWilliams}, we deduce that there exists an $(\Fc_t)$-predictable process $\hat u$ with values in $U$ such that $ u=\hat u$, $d\P\otimes dt$-a.e., so, in particular, $J^\beta(x, u)=J^\beta(x,\hat u)$. By Lemma 2.17-b) in \cite{JacodShiryaev}, it follows that there exists an $(\tilde\Fc_t)$-predictable process $\tilde u$ with values in $U$ which is indistinguishable from $\hat u$, so that $J^\beta(x,\hat u)=J^\beta(x,\tilde u)$. In addition, since $\tilde u$ is in particular $(\tilde\Fc_t)$-progressively measurable, for every $\omega^2\in\bar\Omega^2$ we have that the process $\tilde u^{\omega^2}$ on $(\bar\Omega^1,\bar\Fc^1,\bar\P^1)$, given by $\tilde u_t^{\omega^2}(\omega^1):=\tilde u_t(\omega^1,\omega^2)$, is $(\bar\Fc_t^1)$-progressively measurable. In other words, $\tilde u^{\omega^2}\in\bar\Uc$ for every $\omega^2\in\bar\Omega^2$.

Consider now, for every $\omega^2\in\bar\Omega^2$, the process $(\bar X^{x,\tilde u^{\omega^2}})_{t\geq0}$ solving the following controlled equation:
\[
\bar X_t^{x,\tilde u^{\omega^2}} \ = \ e^{tA} x  + \int_0^t e^{(t-s)A} F\big(\bar X_s^{x,\tilde u^{\omega^2}},\tilde u_s^{\omega^2}\big) \, ds +  \int_0^t e^{(t-s)A}G \, d\bar W_s^1, \qquad \text{for all }t\geq0,\,\bar\P^1\text{-a.s.}
\]
On the other hand, we recall that the process $(X^{x,\tilde u})_{t\geq0}$ solves the controlled equation
\[
X_t^{x,\tilde u} \ = \ e^{tA} x  + \int_0^t e^{(t-s)A} F(X_s^{x,\tilde u},\tilde u_s) \, ds +  \int_0^t e^{(t-s)A}G \, dW_s^1, \qquad \text{for all }t\geq0,\,\P\text{-a.s.}
\]
So, in particular, there exists a $\P$-null set $N\subset\Omega$ such that the above equality holds, for all $t\geq0$ and for every $\omega\notin N$. Therefore, there exists a $\bar\P^2$-null set $\bar N^2\in\bar\Fc^2$ such that, for every $\omega^2\notin\bar N^2$,
\[
X_t^{x,\tilde u}(\cdot,\omega^2) \ = \ e^{tA} x  + \int_0^t e^{(t-s)A} F(X_s^{x,\tilde u}(\cdot,\omega^2),\tilde u_s(\cdot,\omega^2)) \, ds +  \int_0^t e^{(t-s)A}G \, dW_s^1, \; \text{for all }t\geq0,\,\bar\P^1\text{-a.s.},
\]
which can be rewritten in terms of $\tilde u^{\omega^2}$ as
\[
X_t^{x,\tilde u}(\cdot,\omega^2) \ = \ e^{tA} x  + \int_0^t e^{(t-s)A} F(X_s^{x,\tilde u}(\cdot,\omega^2),\tilde u_s^{\omega^2}) \, ds +  \int_0^t e^{(t-s)A}G \, dW_s^1, \quad \text{for all }t\geq0,\,\bar\P^1\text{-a.s.}.
\]
Then, we see that, for every $\omega^2\notin\bar N^2$, the two processes $(\bar X^{x,\tilde u^{\omega^2}})_{t\geq0}$ and $(X^{x,\tilde u}(\cdot,\omega^2))_{t\geq0}$ solve the same equation. By pathwise uniqueness, it follows that, for every $\omega^2\notin\bar N^2$, $(\bar X^{x,\tilde u^{\omega^2}})_{t\geq0}$ and $(X^{x,\tilde u}(\cdot,\omega^2))_{t\geq0}$ are indistinguishable. An application of Fubini's Theorem yields
\[
J^\beta(x,\tilde u) \ = \
\int_{\bar\Omega^2} \bar\E^1\bigg[\int_0^\infty e^{-\beta t} \ell\big(\bar X_t^{x,\tilde u^{\omega^2}},\tilde u_t^{\omega^2}\big)\,dt
\bigg] \, \bar\P^2(d\omega^2) \ = \ \bar\E^2\big[\bar J^\beta\big(x,\tilde u^{\omega^2}\big)\big] \ \geq \ \bar v^\beta(x).
\]
Recalling that $J^\beta(x, u)=J^\beta(x,\tilde u)$, the claim follows taking the infimum over all $ u\in\Uc$.
\ep

\section{Randomized optimal control problems}
\label{S:RandInfHor}

\setcounter{equation}{0} \setcounter{Assumption}{0}
\setcounter{Theorem}{0} \setcounter{Proposition}{0}
\setcounter{Corollary}{0} \setcounter{Lemma}{0}
\setcounter{Definition}{0} \setcounter{Remark}{0}

In the present section we formulate the randomized versions (see \cite{KharroubiPham}) of both the infinite and the finite horizon stochastic optimal control problems introduced in the previous Section \ref{S:InfHor}.

We consider the same probabilistic setting as in subsection \ref{SubS:InfHorProductSpace}. In particular, we adopt the same notations: $(\Omega,\Fc,\P)$, $W^1$, $W^2$, $(\Fc_t)$, $\Uc$. Progressive measurability of processes will always be intended with respect to the filtration $(\Fc_t)$.

By $\Acn$ we denote the family of progressively measurable   processes $\alpha$ with values in $U$  such that
$|\alpha|\leq n$, $\mathbb{P}\otimes dt$-almost surely.
Moreover $\Ac:=\cup_{n\in \mathbb{N}} \Acn$ is the set of 
 progressively measurable and essentially  bounded processes with values in $U$.

\paragraph{State process.} Given $(x,a)\in H\times U$ and $\alpha\in\Ac$, we  consider the  system of controlled stochastic differential equations:
\begin{equation}\label{System}
\begin{cases}
d\X_t\ &= \ A\X_t dt + F(\X_t, \I_t)dt +\,GdW_t^1, \hspace{8.7mm} \X_0 \ = \ x, \\
d \I_t\ &= \ R\alpha_t dt + RdW_t^2,  \hspace{33.8mm} \I_0 \ = \ a. \\
\end{cases}
\end{equation}
On $F$ and $G$ we impose the same assumptions as in Section \ref{S:InfHor}, while on  $R$ we impose the following:

\begin{itemize}
\item[({\bf A.8})] $R\colon U\rightarrow U$ is a trace class injective linear operator with dense image.
\end{itemize}

\begin{Proposition}\label{Eq-stato.randomizzata}
Assume {\bf(A.1)}--{\bf(A.5)} and {\bf(A.8)}. Then, for any $(x,a)\in H\times U$ and $\alpha\in\Ac$, there exists a unique (up to indistinguishability) pair of processes $\X^{x,a,\alpha}=(\X_t^{x,a,\alpha})_{t\geq0}$ and $\I^{a,\alpha}=(\I_t^{a,\alpha})_{t\geq0}$ (the process $\I^{a,\alpha}$ is independent of $x$) such that:
\begin{itemize}
\item $\I^{a,\alpha}$ is given by:
\begin{equation}\label{def-di-I}
\I_t^{a,\alpha} = a  + \int_0^t R\alpha_s \, ds +  RW_t^2, \qquad \text{for all }t\geq0,\,\P\text{-a.s.}
\end{equation}
thus satisfies the second equation in \eqref{System} and
 belongs to $L^{p,loc}_\Pc (\Omega;C([0,+\infty[;U))$ for all $p \geq 1$;
\item $\X^{x,a,\alpha}$ belongs to $L^{p,loc}_\Pc (\Omega;C([0,+\infty[;H))$ for all $p \geq 1$
and  is a mild solution of the first equation in \eqref{System}, that is:
\[
\X_t^{x,a,\alpha} = e^{tA} x  + \int_0^t e^{(t-s)A} F(\X_s^{x,a,\alpha},\I_s^{a,\alpha}) \, ds +  \int_0^t e^{(t-s)A}G \, dW_s^1, \qquad \text{for all }t\geq0,\,\P\text{-a.s.}
\]
\end{itemize}
Moreover for every $0<\rho<\frac{1}{2}-\gamma$ (with $\gamma$ as in Assumption {\bf (A.3)}), and for any $T>0$, the following estimate holds:
\begin{equation}\label{stima_compattezza}
\sup_{t \in [0,T]} t^{\rho}\,\E || \X_t^{x,a,\alpha} ||_{\mathcal{D}((\delta I-A)^{\rho})} \ \leq \ c, 
\end{equation}
for some positive constant $c$, depending only on $T$, $x$, $\rho$, and on the constants introduced in Assumptions {\bf(A.1)}--{\bf(A.5)}, but independent of $a$ and $\alpha$.

\end{Proposition}
\Dim
This is quite a classical result, see for instance \cite{DpZ1} and the randomization framework has nothing special here (we just formulate the result in the case in which we need it). For the sake of completeness, we report the proof of estimate \eqref{stima_compattezza} . We have
\begin{align*}
\E\|\X_t\|_{ \mathcal{D}(A^{\rho})} &\leq \E\|e^{tA} x\| _{ \mathcal{D}(A^{\rho})} +  \E \left\|\int_0^t  e^{(t-s)A}  F(\X_s,\I_s) \, ds \right\| _{ \mathcal{D}((\delta I-A)^{\rho})}\!\!\!\!\!\!\!\!+ \E\left\|\int_0^t e^{(t-s)A}  G \, dW_s^1\right\| _{ \mathcal{D}((\delta I-A)^{\rho})} \\
&\leq \frac{|x|_H}{t^\rho} +  k \int_0^t (t-s)^{-\rho} \E|F(\X_s,\I_s)|_H \, dt +  \bigg( \E \bigg\| \int_0^t e^{(r-s)A}G dW_s^1 \bigg\| _{ \mathcal{D}((\delta I-A)^{\rho})}^2\bigg)^{\frac{1}{2}}
 \\
&\leq \frac{|x|_H}{t^\rho} +  k C_F t^{1-\rho}(1+ \E \sup_{ t \in [0,T]}|\X_t|_H) +  L\Big(  \E \int_0^t (t-r)^{-2\rho -2\gamma}\, dr \Big)^{1/2}\\
&\leq \frac{|x|_H}{t^\rho} +  k C t^{1-\rho}  + L  t ^{\frac{1}{2}-\rho -\gamma}.
\end{align*}
Thus, for every $\rho>0$ such that $ \rho +\gamma < \frac{1}{2}$, we deduce estimate \eqref{stima_compattezza}.
\finedim

\begin{Remark}\label{R:J^R=J}
Notice that $\X^{x,a,\alpha}=X^{x,\I^{a,\alpha}}$, with $X^{x,u}$ defined in Proposition \ref{RegStatoU}. Indeed, for every $a \in U$ and $\alpha\in \mathcal{A}$,  the process $\I^{a,\alpha}$ belongs to  $\mathcal{U}$. Thus the analogue of estimates \eqref{stimadisstato-oriz-fin-prod} and \eqref{stimadisstato-oriz-infin-prod} hold for $\X^{x,a,\alpha}$ (uniformly with respect to $a$ and $\alpha$).
\end{Remark}

\medskip

\noindent Once more we define, in this new setting, the finite and infinite horizon costs as well as the corresponding value functions.
 
\paragraph{Infinite horizon control problem.} For every $\beta>0$, and for any $(x,a)\in H\times U$, $\alpha\in\Ac$, the infinite horizon cost functional is
\[
J^{\beta,\Rc}(x,a,\alpha) \ := \ \E\bigg[\int_0^\infty e^{-\beta s}\,\ell(\X_s^{x,a,\alpha},\I_s^{a,\alpha})\,dt\bigg],
\]
with corresponding value function
\[
v^{\beta,\Rc}(x,a) \ := \ \inf_{\alpha\in\Ac}J^{\beta,\Rc}(x,a,\alpha).
\]

\paragraph{Finite horizon control problem.} For every $T>0$, and for any $(x,a)\in H\times U$, $\alpha\in\Ac$, the finite horizon cost cost is
\[
J^{\beta,T,\Rc}(x,a,\alpha) \ := \ \E\bigg[\int_0^T e^{-\beta s} \ell(\X_{s}^{x,a,\alpha},\I_{s}^{a,\alpha})\,ds+e^{-\beta T}\phi(\X_{T}^{x,a,\alpha})\bigg],
\]
with corresponding the value function
\[
v^{\beta,T,\Rc}(x,a) \ := \ \inf_{\alpha\in\Ac}J^{\beta,T,\Rc}(x,a,\alpha).
\]
Next statement entitles us to study the (asymptotic) behavior of $v^{\beta,\Rc}$ and $v^{T,\Rc}$ instead of $\bar v^{\beta}$ and $\bar v^{T}$ (or of $ v^{\beta}$ and $ v^{T}$ ). Moreover it implies that $v^{\beta,\Rc}$ and $v^{T,\Rc}$ do not depend on their last argument. 
\begin{Proposition}\label{P:v=V^R}
Suppose that Assumptions {\bf(A.1)}--{\bf(A.8)} hold. Then, we have (recalling Proposition \ref{UguaglianzaValueFunction1}):
\begin{itemize}
\item[\textup{(i)}] For every $\beta>0$, $\bar v^\beta(x)=v^\beta(x)=v^{\beta,\Rc}(x,a)$, for all $(x,a)\in H\times U$. In par\-ti\-cu\-lar, the function $v^{\beta,\Rc}$ is independent of its second argument.
\item[\textup{(ii)}] For every $T>0$ and $\beta\geq 0$, $\bar v^{\beta,T}(x)=v^{\beta,T}(x)=v^{\beta,T,\Rc}(x,a)$, for all $(x,a)\in H\times U$. In par\-ti\-cu\-lar, the function $v^{\beta,T,\Rc}$ is independent of its second argument.
\end{itemize}
\end{Proposition}

\Dim
We only report the proof of the first statement, as item (ii) can be proved in an analogous way. We split the proof of (i) into some steps.

\vspace{2mm}
 
\noindent{\bf Step 1.} We are going to prove that,  for every fixed $T>0$, $x \in H$ and  $u \in \mathcal{U}$,
\begin{equation}\label{step1}
 \text{If} \ u^n\xrightarrow{\mathbb{P}} u  \ \text{ in }  \Omega \times [0,T], \ \text{ then }  X^{x,u^n} \rightarrow  X^{x,u} \text{ in } L^2_{\mathcal{P}}(\Omega \times [0,T];H)
\end{equation}
We set $\bar{X}^n= X^{x,u^n}-  X^{x,u}$, thus 
\[  \bar{X}^n _t= \int_0^t e^{(t-s)A} [F(X_s^{x,u^n},u^n_s )- F(X_s^{x,u}, u_s)] \, ds , \qquad \text{for all }t\in [0,T],\,\P\text{-a.s.}
 \]
and  there is a positive constant $C$ (that depends only on  the Lipschitz constant of $F$ and on $T$) such that
\[
| \bar{X}^n _t|^2 \leq  C \left[\int_0^t  |\bar{X}^n_s|^2 \, ds +  \int_0^t |F(X_s^{x, u},u^n_s )- F(X_s^{x,u},u_s)| ^2\, ds \right].
\]
Thus:
\[
\E \sup_{t \in [0,r]}| \bar{X}^n _t|^2 \leq  C \int_0^r \E \sup_{\sigma \in [0,s]}  |\bar{X}^n_\sigma|^2 \, ds + C \,  \E  \int_0^T |F(X_s^{x, u},u^n_s )- F(X_s^{x,u},u_s)| ^2\, ds.
\]
Hence, by Gronwall's lemma
\[
\E \sup_{t \in [0,r]}| \bar{X}^n _t| ^2\leq   C e^{CT} \, \E  \int_0^T |F(X_s^{x, u},u^n_s)- F(X_s^{x,u},u_s)| ^2\, ds.
\]
Notice that by {\bf{(A.4)}}  we have 
\[    |F(X_s^{x, u},u^n_s )- F(X_s^{x,u},u_s)| \leq 2 C_F (1 + |X_s^{x,u}|)\]
Henceforth, thanks to Proposition \ref{RegStatoU}, we can apply Lebesgue's dominated convergence Theorem  to derive \eqref{step1}. Observe that, actually, we have proved the following:
\[
\lim_{n\to +\infty }\E \sup_{t \in [0,r]}| X^{x,u^n}_t-  X^{x,u}_t|^2  \to 0.
\]

\vspace{2mm}
 
\noindent{\bf  Step 2.} Fix  $x \in H$, $u \in \mathcal{U}$. We are going to show that: 
 \begin{equation}\label{step2}
  \text{If, $\forall\, T>0$,} \ u^n\xrightarrow{\mathbb{P}\otimes dt} u  \ \text{ in }  \Omega \times [0,T], \ \text{ then } \lim_{n \to \infty} J^\beta(x,u^n )= J^\beta(x,u).
\end{equation}
Thanks to the presence of the discount term $e^{-s\beta}$ and the boundedness of $\ell$, see ${\bf(A.6)}$,  it is enough to prove that 
\[  \lim_{n \to \infty} \E\int_0^T \ell(X^{x,u^n}_s,u^n_s) \, ds =\E  \int_0^T\ell(X^{x,u}_s,u_s) \, ds. \]
Indeed, for every $T>0$,
\[  |J^\beta({x,u^n})-  J^\beta({x,u}) | \leq   \left|\E\int_0^T [\ell(X^{x,u^n}_s,u^n_s)- \ell(X^{x,u}_s,u_s)] \, ds\right| + 2M_\ell \int_T^{+\infty} e^{-\beta s} \, ds. \]
From  {\bf Step 1} and the boundedness of $\ell$, we obtain
\[  \lim_{n \to \infty} \E\int_0^T \ell(X^{x, u^n}_s,u^{n}_s) \, ds =\E  \int_0^T\ell(X^{x,u}_s,u_s) \, ds. \]
By the previous considerations we deduce the validity of \eqref{step2}. 

\vspace{2mm}
 
\noindent{\bf Step 3.} Let $\Uc_{\textup{bdd}}$ denote the subset of $\Uc$ of all uniformly bounded processes. By the previous step, we deduce the following equality:
\begin{equation}\label{Ubdd}
\inf_{u\in\Uc}J^\beta(x,u) \ = \ \inf_{u\in\Uc_{\textup{bdd}}}J^\beta(x,u), \qquad \text{for every }x\in H,
\end{equation}
As a matter of fact, given $u\in\Uc$, it is enough to apply \eqref{step2} to the sequence $ u^n = I_{B(0,n)}(|u|)u$, $n \in \mathbb{N}$.

\vspace{2mm}
 
\noindent{\bf Step 4.} Fix $a\in U$. Given $u\in\Uc_{\textup{bdd}}$, we claim that there exists a sequence $(\alpha^n)_n\subset\Ac$ such that (see \eqref{def-di-I} for the definition of $\I^{a,\alpha^n}$)
\[
\I^{a,\alpha^n} \to u \qquad \text { in } L^2( [0,T] \times \Omega;U) \qquad \forall\,T >0.
\]
The above follows if we prove that the affine set
$$\{\I^{a,\alpha}\colon \alpha\in\Ac\}$$
is dense in $L^2_{\mathcal{P}}( [0,T] \times \Omega; U)$. Since $\I^{a,\alpha}:=\hat{\I}^{\alpha}+a + RW^2$ with $\displaystyle \hat{\I}^{\alpha}_t:=  \int_0^t R\, \alpha _s \, ds$, it is enough to prove that the linear subspace
$\displaystyle \{\hat{\I}^{\alpha}\colon \alpha\in\Ac\}$
is dense in $L^2_{\mathcal{P}}( [0,T] \times \Omega; U)$.

Now, consider the linear space generated by the set of functions 
\[
\{ \eta {I}_{[t_0,T)}: t_0 \in [0,T), \; \eta: \Omega \to U \; \mathcal{F}_{t_0}\text{-meas. and bounded}\}.
\]
As it is well-known such linear space is dense in $L^2([0,T] \times \Omega; U)$ (notice that it coincides with the linear space generated by the set of functions of the form: $\eta {I}_{[t_0,t_1)}$, with $0\leq t_0<t_1\leq T$, $ \eta  $ $ \mathcal{F}_{t_0}$-measurable and bounded).

Thus it is enough to prove that, for all $t_0\in [0,T)$ and every $\mathcal{F}_{t_0}$-measurable and bounded  $\eta$, there exists a sequence 
 $(\alpha^n )_n\subset\Ac$,  such that:
$$\mathbb{E}\int_0^T|\eta {I}_{[t_0,T)}(t)-\hat{\I}^{\alpha^n}_t|^2 dt\rightarrow 0.$$

Now, let $(U_m)_m$ be a sequence of finite dimensional subspaces of $U$ such that $U_{m}\subset U_{m+1}$ and ${\cup_{m=1}^{\infty} U_m}=U$ and define $E_m=R U_m$. Clearly $E_{m}\subset E_{m+1}$, ${\cup_{m=1}^{\infty} E_m}$ is dense in $U$ and $R$ is invertible from  $U_m$ to $E_m$ (thus bounded with bounded inverse).

We can always assume that $\eta$ takes values in $E_m$ for a suitable $m$, possibly approximating $\eta$ by its projection on $E_m$.

Now take $\alpha^n_s =  n {I}_{[t_0,t_0+\frac{1}{n}]} (s) R^{-1}\eta$, $ n \in \mathbb{N}$, so that  for every $n$ we have that $\alpha_n\in\Ac$. Moreover
$\hat{\I}^{\alpha^n}=n(1-(t-t_0))\eta {I}_{[t_0,t_0+1/n[}+\eta {I}_{[t_0+1/n,T]}$ and consequently
$$\mathbb{E}\int_0^T|\eta {I}_{[t_0,T)}(t)-\hat{\I}^{\alpha^n}_t|^2 dt= 
\frac{1}{n}\mathbb{E}|\eta|^2.$$

\vspace{2mm}
 
\noindent{\bf Conclusion.} Recalling from Remark \ref{R:J^R=J} that $J^{\beta,\Rc}(x,a,\alpha) = J^{\beta}(x,\I^{a,\alpha})$, we immediately see that $v^{\beta}(x)\leq v^{\beta,\Rc}(x,a)$.

On the other hand, given $u\in\Uc_{\textup{bdd}}$, by  \textbf{Step 4} and a standard diagonal argument it is possible to construct a sequence $(\alpha^n)_n\in\Ac$  such that 
$\forall\, T>0$, $\I^{a,\alpha^n}\xrightarrow{\mathbb{P}\otimes dt} u$  in $  \Omega \times [0,T]$. Thus, by \textbf{Step 2},
$$v^{\beta,\Rc}(x,a)\leq \lim_{n\rightarrow \infty} J^{\beta,\Rc}(x,a,\alpha^n) = \lim_{n\rightarrow \infty} J^{\beta }(x,\I^{a,\alpha^n})=
J^{\beta}(x,u).$$ 
Hence, using also \eqref{Ubdd}, we see that the reverse inequality  $v^{\beta}(x)\geq v^{\beta,\Rc}(x,a)$ holds as well. \finedim

\section{BSDE representation of the value functions}
\label{S:BSDE}

In the present section we obtain a non-linear Feynman-Kac formula for the value functions of both the infinite and the finite horizon stochastic optimal control problems.
This representation will be the essential tool to study the asymptotic behavior of the value functions.

\subsection{Elliptic BSDE: infinite horizon optimal control problem}

We denote by $\X^{x,a}$  and $\I^a$ the solution to \eqref{System} corresponding to $\alpha \equiv 0$.
For any $\beta >0$ and $n \in \mathbb{N}$, we consider the following  standard BSDE with infinite terminal time and generator being Lipschitz with respect to the martingale variable $(Z,\Gamma)$ and strongly dissipative with respect to $Y$:

\begin{align} \label{EllitticaRandPen}
Y_t^{x,a,\beta,n} &=  Y_T^{x,a,\beta,n}  -\beta \int_t^T Y_s^{x,a,\beta,n} \, ds + \int_t^T \ell(\X^{x,a}_s, \I^{a}_s)\,ds - n \int_t^T |\Gamma^{x,a,\beta,n}_s| \, ds \nonumber\\ &\quad-\int_t^T Z_s^{x,a,\beta,n} \, d W^1_s 
-\int_t^T \Gamma^{x,a,\beta,n}_s \, d W^2_s,  \qquad 0 \leq t \leq T <+\infty.
\end{align}
\begin{Proposition}\label{PenEllRand}
Suppose that Assumptions {\bf (A.1)}--{\bf(A.8)} hold. Then
\begin{enumerate}
\item [(i)] For every $(x,a) \in H \times U$, there exists a unique solution $ (Y^{x,a,\beta,n}, \Gamma^{x,a,\beta,n}, Z^{x,a,\beta,n})$ of the BSDE \eqref{EllitticaRandPen}
with $Y^{x,a,\beta,n}$ bounded, continuous and progressively measurable, $ Z ^{x,a,\beta,n}$ belonging to $  L^{2, loc}_{\mathcal{P}}(\Omega\times [0,+\infty[;\Xi^*)$ and $\Gamma ^{x,a,\beta,n}$ belonging to $ L^{2, loc}_{\mathcal{P}}(\Omega\times [0,+\infty[;U^*)$.
\item[(ii)] The following bounds hold (uniformly with respect to $n$):
\begin{equation}\label{stimaYpenunif}
| Y^{x,a,\beta,n}_t | \leq \frac{M_\ell}{\beta},
\end{equation}
\begin{equation}\label{stimaZpenunif}
| Z^{x,a,\beta,n}_t | \leq \frac{L_\ell |G|}{\mu},
\end{equation}
\begin{equation}\label{stimaUpenunif}
 \E \int_0^{+\infty}e^{-2\beta s} |\Gamma_s^{x,a,\beta,n}|^2 \, ds  < \infty.
\end{equation}
\item[(iii)] For  every $(x,a)\in H\times U$,  if we define the value function for the infinite horizon  problem, in the randomized framework and  with bounded set of controls, namely
\begin{equation}\label{defvbetan}
 v^{\beta,n,\Rc}(x,a):= \ \inf_{\alpha\in\Acn}J^{\beta,\Rc}(x,a,\alpha),
\end{equation}
then, for all $t\geq 0$:
\begin{equation}\label{rapprsolbetan}
Y_t^{x,a,\beta,n} = v^{\beta,n,\Rc}(\X^{x,a} _t, \I^{x,a} _t), \qquad \P\text{- a.s.}
\end{equation}
In particular 
\begin{equation}\label{Yfunzionevaloreconn}
Y^{x,a,\beta,n}_0 \ = \ v^{\beta,n,\Rc}(x,a) \ = \ \inf_{\alpha\in\Acn}J^{\beta,\Rc}(x,a,\alpha).
\end{equation}
\end{enumerate}
\end{Proposition}
\Dim
Equation \eqref{EllitticaRandPen}  fulfills the standard assumptions in \cite{HuTess}, Lemma 2.1,   thus we already know that  $(i)$, as well as estimates \eqref{stimaYpenunif} and \eqref{stimaUpenunif} in $(iii)$, hold true.

It remains to prove the uniform estimate \eqref{stimaZpenunif}. To do that we need to introduce finite horizon approximations  of \eqref{EllitticaRandPen} and then to smooth up its coefficients.

Denote by  $(Y^{x,M}, \Gamma^{x,M}, Z^{x,M})$ the  solution to the finite horizon BSDE on $[0,M]$, with $M\in\N$: 
\begin{align*} 
Y_t^{x,M} =  & -\beta \int_t^M Y_s^{x,M} \, ds + \int_t^M \ell(\mathfrak{X}^{x,a}_s, \I^{a}_s)\,ds - n \int_t^M |\Gamma^{x,M}_s| \, ds \nonumber\\ &-\int_t^M Z_s^{x,M} \, d W^1_s 
-\int_t^T \Gamma^{x,M}_s \, d W^2_s,  \forall\,t\in[0,M],
\end{align*}
where we have omitted some parameters in the notation to keep it readable. We have (see again \cite{HuTess} Lemma 2.1, ):
$$\mathbb{E}\int_0^T |Z^{x,M}_t-Z_t^{x,a,\beta,n}|^2
 dt \rightarrow 0 \;\; \hbox{ as } M\rightarrow \infty.
$$
Thus it is enough to prove \eqref{stimaZpenunif} for $Z^{x,M}$. To this end, we construct sequences $(F_k(\cdot))_{k\in \mathbb{N}}$ and $(\ell_k(\cdot,\cdot))_{k\in \mathbb{N}}$ of Gateaux differentiable functions punctually converging to $F$ and $\ell$  such that assumptions \textbf{(A.4)}, \textbf{(A.5)} and  \textbf{(A.6)} hold with the same constants {(see \cite{BismutElworthy} for such a construction)}.

We introduce the notation $\rho_k(u)=\sqrt{|u|^2+k^{-1}}$.  Consider the regularized forward-backward system parametrized with final time $M$ (we omit parameters in the equations to simplify notation):
\begin{equation}\begin{cases}\label{systemF_k}d\X_t\ &= \ A\X_t dt + F_k(\mathfrak{X}_t, \I_t)dt +\,GdW_t^1, \hspace{37.4mm} \X_0=x,\; t\in [0,M], \\
d \I_t\ &= RdW_t^2,  \hspace{64.7mm} \I_0=a,\; t\in [0,M],\\
dY_t\ &= \ -\beta Y_tdt- \ell_k(\mathfrak{X}_{t},\I_{t})dt-n\rho_k(\Gamma_t)dt-Z_tdW^1_t-\Gamma_tdW^2_t,  \hspace{9.7mm} Y_{M}=0.
\end{cases}\end{equation}
Denote by $(\X^{x,k}, \I, Y^{x,k,M},Z^{x,k,M},\Gamma^{x,k,M})$ the solution to the above forward-backward system  (for readability sake we do not report parameter $a$). 
 Proceeding as in \cite{BismutElworthy}, it is easy to verify by parameter depending contraction theorem that, as $k\rightarrow \infty$,
$$ 
\mathbb{E}\int_0^M|Z^{x,k,M}_t-Z^{x,M}_t|^2 dt\rightarrow 0.$$
Thus, once more,  it is enough to prove  \eqref{stimaZpenunif} for $Z^{x,k,M}$. 

Set $\nu^M(\tau,x,a)= Y^{x,k,M-\tau}_{0}$. Then, recalling that all the coefficients in \eqref{systemF_k} are differentiable we have the following identifications (see \cite{FuTessAOP2002}) 
$$Y^{x,k,M}_t=\nu^M(t,\mathfrak{X}^{x,k}_t,\I_t),\qquad Z^{x,k,M}_t=\nabla_x\nu^M(t,\mathfrak{X}^{x,k}_t,\I_t)\,G. $$
Thus, to obtain \eqref{stimaZpenunif} it is enough to prove that, for all $M >0$
\[
|\nabla_x Y^{x,k,M}_{\tau}|\leq \frac{L_{\ell}}{\beta}.
\]
Differentiating  \eqref{systemF_k} with respect to $x$ in the direction $\xi$, 
 we get
$$\begin{cases}d\nabla^{\xi}_x \mathfrak{X}_t\ &= \ A\nabla^{\xi}_x \mathfrak{X}_t dt + \nabla_x F_k(\mathfrak{X}_t, \I_t)\nabla^{\xi}_x \mathfrak{X}_tdt, \\
 \nabla^{\xi}_x \mathfrak{X}_\tau \ &= \ \xi,\\
-d\nabla^{\xi}_x Y_t\ &= \ -\beta \nabla^{\xi}_x Y_tdt- \nabla_x \ell_k(\mathfrak{X}_t,\I_t)\nabla^{\xi}_x \mathfrak{X}_tdt -n\nabla^{\xi}_x \rho_k(U_t) \nabla^{\xi}_x U_tdt-\nabla^{\xi}_x Z_tdW^1_t-\nabla^{\xi}_x \Gamma_tdW^2_t,\\ \nabla^{\xi}_x Y_{\ell}\ &= \ 0
. \end{cases}
$$
Exploiting the dissipativity of $A$ and $F_k$ (namely \textbf{(A.1)} and \textbf{(A.4)}), 
 (see, for instance \cite{HuTess}) we get $|\nabla_x^{\xi} \X_t|\leq e^{-\mu (t-\tau)} |\xi|$
and,  again by a standard Girsanov argument, we deduce that, with respect to a  suitable probability $\tilde{\mathbb{P}}$,
$$|\nabla_x^{\xi} Y_{\tau}|= \left|\tilde{\mathbb{E}}\int_\tau^{\ell} e^{-\beta (t-\tau)}\nabla_x \ell_k(\mathfrak{X}_t,\I_t)\nabla^{\xi}_x \mathfrak{X}_t dt\right| \leq \frac{L_\ell}{\mu+\beta}. $$
Thus estimate \eqref{stimaZpenunif} holds.

To prove \eqref{Yfunzionevaloreconn} we have to come back to the results in \cite{HuTess} paying some attention to the fact that here we want to work with a strong formulation of the control problem (that is with fixed probability space and noise).
\newline If we take into account that the Hamiltonian function is, in this case, given by:
$$\inf_{|\alpha|\leq n}\{\ell(x,u)+u \alpha\}=
\ell(x,u)-n|u|,$$ the fundamental relation (see relation (5.12) in \cite{HuTess})  provides, computing expectation and letting $T\rightarrow \infty$,
$$J^{\beta,\Rc}(x,a,\alpha)=Y^{x,a,\beta,n}_0+\mathbb{E}\int_0^{\infty} e^{-\beta s}\left [ \Gamma ^{x,a,\beta,n}_s \alpha_s -n |\Gamma ^{x,a,\beta,n}_s |\right] ds$$
and the claim follows letting $\alpha_s=n \Gamma ^{x,a,\beta,n}_s |\Gamma ^{x,a,\beta,n}_s |^{-1}$.

\noindent Finally, \eqref{rapprsolbetan} is a straightforward  consequence of \eqref{Yfunzionevaloreconn} and of the Markovianity of the problem.
\finedim

Taking into account the definition of $v^{\beta,n,\mathcal{R}} (x,a)$ and identification \eqref{defvbetan}, Propositions  \ref{UguaglianzaValueFunction1} and \ref{PenEllRand} yield:

\begin{Proposition} \label{proprietan}
If Assumptions {\bf (A.1)}--{\bf(A.8)}  are verified, then the following properties hold for every $x,x'\in H$, $a\in A$:
\begin{align}
 v^{\beta,n,\mathcal{R}}(x,a) \ \leq \ \frac{M_\ell}{\beta}, \quad & 
|  v^{\beta,n,\mathcal{R}}(x,a)-  v^{\beta,n,\mathcal{R}}(x',a)| \leq \ L_\ell \, |x-x'|,\label{unifLipsn} \\
v^{\beta,n,\Rc}(x,a) \ & \downarrow  \ \inf_{\alpha\in\Ac }J^{\beta,\Rc}(x,a,\alpha) \ = \ v^{\beta, \mathcal{R}}(x) = \ v^{\beta}(x), \label{decreasing} \\
|v^{\beta}(x)| \leq \ \frac{M_\ell}{\beta}, \ & \quad 
|  v^{\beta}(x)-  v^{\beta}(x')| \leq \ L_\ell \, |x-x'|,\label{unifLipsn-ell}
\end{align}
\end{Proposition}

\vspace{3mm}

\noindent Our aim is to characterize $v^{\beta,\Rc}$ in terms of the maximal  solution to the following  non standard `constrained' infinite horizon BSDE involving a  reflection term: 
\begin{equation}\label{EllitticaRand}
Y_t = Y_T  -\beta \int_t^T Y_s \, ds + \int_t^T \ell(\X^{x,a}_s, \I^{a}_s)\, ds + K_t - K_T- \int_t^T Z_s  \, d W^1_s, \quad \forall\, 0 \,  \leq t \leq T,\; \ \mathbb{P}\text{-a.s.}
\end{equation}
\begin{Definition}\label{SolEllRand}
A solution to \eqref{EllitticaRand}  is a triple  $(Y, Z, K)$ such that $Y$ is a progressively measurable bounded process with continuous trajectories, $ Z \in  L^{2, loc}_{\mathcal{P}}(\Omega\times [0,+\infty[;\Xi^*)$,  $ K \in \mathcal{K}^{2,loc}$ and
 \eqref{EllitticaRand} holds.
\end{Definition}
\begin{Proposition}\label{propEllitica}
Assume  {\bf(A.1)}--{\bf(A.8)}. Then, for every $\beta  >0$ and $(x,a) \in H\times U$, there exists a triple 
$(Y^{x,a,\beta}, Z^{x,a,\beta}, K^{x,a,\beta})$ of progressively measurable processes such that:

\begin{enumerate}
\item $Y^{x,a,\beta}$ is the decreasing limit of $(Y^{x,a,\beta,n})_{n}$, see \eqref{EllitticaRandPen},
moreover the following holds, for all $t\geq0$,
\begin{equation}\label{CondEll}
 Y^{x,a,\beta}_t = v^\beta (\X^{x,a}_t),\quad |Y^{x,a,\beta}_t | \leq \frac{M_\ell}{\beta}, \qquad \mathbb{P}\text{-a.s.}.
\end{equation}

\item The following estimate holds, for all $t\geq0$,
 \begin{equation}\label{CondZEll} 
|Z^{x,a,\beta}_t | \leq \frac{L_\ell}{\mu}, \qquad  \mathbb{P}\text{-a.s.}
\end{equation} and, for all $T>0$:
\[ 
Z^{x,a,\beta,n} \rightharpoonup Z^{x,a,\beta}  \text { in } L^2_\Pc (\Omega\times [0,T];\Xi ^*).
\]
\item The following convergences take place 
for all $T>0$:
 \[ 
\Gamma^{x,a,\beta,n} \rightharpoonup 0  \text { in } L^2_\Pc (\Omega\times [0,T];U^*),
\]
\[ 
n \int_0^T |\Gamma^{x,a,\beta,n}_s| \, ds \rightharpoonup K^{x,a,\beta}_T \text{in $L^2 (\Omega, \mathcal{F}_t,\mathbb{R})$,} 
\]
\item $({Y}^{x,a,\beta},{Z}^{x,a,\beta},{K}^{x,a,\beta})$ is the maximal solution  of equation \eqref{EllitticaRand} in the sense that if there exists another solution $(\bar{Y},\bar{Z},\bar{K})$ then
$Y_t^{x,a,\beta} \geq \bar{Y}_t^{x,a,\beta}$ for all $t\geq 0$, $\P$-a.s..
\end{enumerate}
\end{Proposition}
\Dim 
The proof is essentially identical to the proof of Proposition 3.3 in \cite{CoFuhPham} and we omit it. 
\finedim

\subsection{Results on the finite horizon case}
\label{sub-finhor}

We report here, for further use, the finite horizon analogue of the results stated in Proposition  \ref{propEllitica}.

\begin{Definition}\label{SolEllRandT}
For every $\beta\geq0$, $T>0$ and $(x,a)\in H\times U$, a solution to the finite horizon BSDE on $[0,T]$
\begin{equation}\label{EllitticaRandT}
Y_t=  \phi(\X^{x,a}_{T})  -\beta \int_t^T Y_s \, ds + \int_t^T \ell(\X^{x,a}_{s}, \I^{a}_{s})\, dr - K_T + K _t- \int_t^T Z_s  \, d W^1_s \ ,
\end{equation}
is a triple $(Y, Z, K)$ in $L^{2}_\Pc (\Omega;{C} ([0,T];\mathbb{R}))\times L^{2}_{\mathcal{P}}(\Omega\times [0,T];\Xi^*)\times \mathcal{K} ^2(0,T)$, satisfying \eqref{EllitticaRandT} for all $t\in [0,T]$, $\mathbb{P}$-a.s..
\end{Definition}

\medskip
Now, for every $n\in\N$, $\beta\geq0$, $T>0$, $(x,a)\in H\times U$, let $(Y^{x,a,\beta,T,n},Z^{x,a,\beta,T,n},\Gamma^{x,a,\beta,T,n})$ $\in$ $L^{2}_{\Pc} (\Omega;{C}( [0,T];\mathbb{R}))\times L^{2}_{\mathcal{P}}(\Omega\times [0,T];\Xi^*)\times  L^{2}_{\mathcal{P}}(\Omega\times [0,T];U^*)$ be the solution  of the standard BSDE on $[0,T]$ (see \cite{FuTessAOP2002}):
\begin{equation}\label{BSDE-fin_hor_n}
Y_t= \phi(\X^{x,a}_T)  -\beta \int_t^T Y_s \, ds + \int_t^T \ell(\X^{x,a}_s, \I^{a}_s)\, ds    - 
 n\int_t^T| \Gamma_s| \, ds -  \int_t^T \Gamma_s \, dW^2_s - \int_t^T Z_s  \, d W^1_s.
\end{equation}
Then, similarly to Proposition \ref{propEllitica}, we have the following result concerning the finite horizon case.

\begin{Proposition}\label{propElliticaT}
Assume  {\bf(A.1)}--{\bf(A.8)}. Then, for every $\beta  \geq0$, $T>0$, and $(x,a) \in H\times U$, there exists 
$(Y^{x,a,\beta,T}, Z^{x,a,\beta,T}, K^{x,a,\beta,T})$ such that:

\begin{enumerate}
\item $(Y^{x,a,\beta,T}, Z^{x,a,\beta,T}, K^{x,a,\beta,T})$ is a solution to \eqref{EllitticaRandT}.
\item $Y^{x,a,\beta,T}$ is the decreasing limit of $(Y^{x,a,\beta,T,n})_{n}$,
and the following representation holds
\[
Y_t^{x,a,\beta,T}=v^{\beta, T-t}(\X_t^{x,a}).
\]
\item $
Z^{x,a,\beta,T,n}$ converges towards $Z^{x,a,\beta,T}$ weakly  in $ L^2_\Pc (\Omega\times [0,T];\Xi^*)$.
\item $
\Gamma^{x,a,\beta,T,n}$ converges towards 0  \text { weakly in } $L^2_\Pc (\Omega\times [t,T];U^*)$.
\item For all $t\in [0,T]$, $n \int_0^t |\Gamma^{x,a,\beta,T,n}_r| \, dr$ converges towards $K^{x,a,\beta,T}_t $ weakly in $L^2 (\Omega, \mathcal{F}_t,\mathbb{R})$.
\item $Y^{x,a,\beta}$ is the maximal solution  of equation \eqref{EllitticaRandT} in the sense that given any other solution $(Y',Z',K')$ of \eqref{EllitticaRandT} then $Y'_t\leq Y_t$, for all $t\in[0,T]$, $\P$-a.s..
\end{enumerate}
\end{Proposition}
\begin{Remark}\label{SolEllRandT0} Concerning the representation at point 2. of Proposition \ref{propElliticaT}, this is a well-known result for the solution to the standard BSDE
\eqref{BSDE-fin_hor_n} and  it is enough to pass to the limit as $n\rightarrow+\infty$.

 Standard estimates on the value function $v^{\beta, T}$ imply that there exists a constant $c_{\beta,T}$, increasing with respect to $T$, such that
$$|Y_t^{x,a,\beta,T}|\leq c_{\beta,T}(1+|\X_t^{x,a}|).$$
\end{Remark}

\noindent To proceed with our arguments we need  to prove that  the maximal solution $(Y^{x,a,\beta}, Z^{x,a,\beta}, K^{x,a,\beta})$ of the infinite horizon BSDE \eqref{EllitticaRand} (see Proposition \ref{propEllitica}), when restricted to an arbitrary compact subset $[0,T]$, coincides with the maximal solution of equation \eqref{EllitticaRandT} with final condition  $\phi$ replaced by $v^{\beta}$ itself, namely of the BSDE:
\begin{equation}\label{EllitticaRandTv}
Y_t=  v^{\beta}(\X^{x,a}_{T})  -\beta \int_t^T Y_s \, ds + \int_t^T \ell(\X^{x,a}_{s}, \I^{a}_{s})\, ds - K_T + K _t- \int_t^T Z_s  \, d W^1_s.
\end{equation}

In the existing literature this result is a straightforward consequence of the characterization of $v^{\beta}$ as the unique viscosity solution of an elliptic HJB equation  (see \cite{CoFuhPham}). Here, where  such tools are not available,  we give a direct proof that avoids the use of PDE techniques.

\begin{Proposition}\label{maximalellitica}
Assume  {\bf(A.1)}--{\bf(A.8)} and fix  $T >0$. Then  $(Y^{x,a,\beta}, Z^{x,a,\beta}, K^{x,a,\beta})$, restricted to $[0,T]$, 
is the maximal solution of the finite horizon BSDE \eqref{EllitticaRandTv}.
\end{Proposition}
\Dim Checking that the restriction of  $(Y^{x,a,\beta}, Z^{x,a,\beta}, K^{x,a,\beta})$ to $[0,T]$ is a solution of equation \eqref{EllitticaRandTv} is straightforward. It remains to prove its maximality.

 To this purpose let 
$(Y^{x,a,\beta,n}, Z^{x,a,\beta,n}, \Gamma^{x,a,\beta,n })$ be the solution of   equation \eqref{EllitticaRandPen}. Moreover  let $(Y^{T}, {Z}^{T}, {K}^{T})$ be a generic solution of the constrained equation \eqref{EllitticaRandTv}. Applying It\^o's formula to $e^{-\beta t} \, (Y^{x,a,\beta,n}_t - Y^{T}_t)$, we obtain
\begin{align*} d [e^{-\beta t} \, (Y^{x,a,\beta,n}_t - Y^{T}_t)]=   &-  e^{-\beta t}( n |\Gamma^{x,a,\beta,n}|)\, dt + e^{-\beta t} d K^{T}_t - e^{-\beta t}(Z^{x,a,\beta,n}_t \! -  Z^{T}_t) \, dW^1_t  \\
& - e^{-\beta t}\Gamma^{x,a,\beta,n} \, d W^2_t.
\end{align*}
Let  $\psi^{\beta,n}$ be defined as follows:
\[
\psi^{\beta,n}_s \ := \ \left\{ \begin{array}{ll}- n  \frac{\Gamma^{x,a,\beta,n}}{|\Gamma^{x,a,\beta,n}|}, & \Gamma^{x,a,\beta,n} \not= 0, \\
0, &   \text{ elsewhere. }\end{array} \right.
\]
Since $\psi^{\beta,n}$ is bounded, its exponential local martingale is indeed  a $\P$-martingale:
\[
L_t^{\beta,n} \ := \ \mathcal{E}\bigg(\int_0^t\psi^{\beta,n}_s \, ds\bigg).
\]
Then, we introduce the new probability $\mathbb{P}^{\beta,n}$ under which $(\tilde {W}^1_t,  \tilde{W}^2_t) =(- \int_0^t \psi^{\beta,n}_s \, ds  + W^1_t,W^2_t)$, is a cylindrical Wiener process (with respect to $(\mathcal{F})$. As a consequence, taking also into account that $K^{T}$ is non-decreasing, we find:
\[
Y^{x,a,\beta,n}_t- Y^{T}_t \geq \E^{\beta,n} \left[\left( e^{-\beta T} ( v^{\beta,n, \mathcal{R}}(\X^{x,a}_T,\I^{a}_T)-  v^\beta(\X^{x,a}_T)) \right)\Big\vert \mathcal{F}_t \right]. 
\]
and the claim follows by \eqref{decreasing}.
\finedim

\section{Asymptotic behavior of the infinite horizon problem and the Ergodic BSDE}
\label{S:Ergodic}
We are now finally  able to introduce the  object of our analysis. For every $(x,a)\in H\times U$, we introduce the following infinite horizon constrained BSDE (ergodic constrained BSDE):
\begin{equation} \label{ErgodicRan}
Y_t \ = \ Y_T + \int_t^T (\ell(\X^{x,a}_s,\I^a_s) -\lambda ) \,ds  - K_T + K_t - \int_t^T Z_s d W^1_s, \qquad 0 \leq t \leq T <+\infty,
\end{equation}
where the real number  $\lambda $ is part of the unknowns.
Such a kind of equation has been already studied in the infinite dimensional setting in \cite {FuhHuTess}, \cite {DebHuTess} under structure condition, while the randomized case  has been addressed in \cite {CoFuhPham},  in the finite dimensional case.

\subsection{Asymptotic behavior of the infinite horizon problem}

As in \cite [Theorem 3.8]{DebHuTess}, we have the following result, which immediately follows from \eqref{unifLipsn-ell}.
\begin{Lemma}\label{Lemmabeta}
Assume {\bf(A.1)}--{\bf(A.8)} and set $\hat{v}^{\beta}(x):={v}^{\beta}(x)-v^\beta(0)$.
Then there exists a sequence $\beta_k \searrow 0$ such that:
 \begin{align}
 \label{defvbar}
\lim_{k\to \infty} \hat{v}^{\beta_k}(x) \ &= \ \hat v(x), \qquad  \forall\,x \in H, \\
 \label{deflambdabar}
\lim_{k\to \infty}\beta_k \hat v^{\beta_k}(0) \ &= \ {\lambda},
\end{align}
for a suitable function $\hat v$ and a suitable real number ${\lambda}$.
\end{Lemma}
\begin{Remark} By \eqref{unifLipsn-ell} and the definition of $\hat v$ we have
\begin{equation}\label{crescita_v} |\hat v^{\beta}(x)|\leq C|x|,\qquad\quad
|\hat v(x)|\leq C|x|,
\end{equation}
where $C$ does not depend on $\beta$.
\end{Remark}

Using the pair $\hat {v}$ and $\lambda$  we can easily construct a solution to equation \eqref{ErgodicRan} (see  \cite{DebHuTess} and \cite{FuhHuTess} for the same argument in the `structure condition' case).

\begin{Theorem}\label{ergodico}
Assume {\bf(A.1)}--{\bf(A.8)}. Let $\hat v$, $\lambda$ and $(\beta_{k})_k$ as in Lemma \ref{Lemmabeta}. Fix $(x,a)\in H\times U$ and  let  $\hat{Y}_t^{x,a}= \hat{v}(X^{x,a}_t)$.\newline There exist  $ \hat{Z}^{x,a}\in L^{2, loc}_{\mathcal{P}}(\Omega\times [0,+\infty[;\Xi ^*)$, $ \hat{K}^{x,a} \in \mathcal{K}^{2,loc}$  and a subsequence $\{ \beta_{k_h} \}$ such that, for all $T>0$,
\begin{align*}
Z^{x,a,\beta_{k_{h}}} &\rightharpoonup \hat{Z}^{x,a}  \text { in } L^2_\Pc (\Omega\times [0,T];\Xi^*), \\
K^{x,a,\beta_{k_{h}}}_T &\rightharpoonup \hat{K}^{x,a} _T \text { in } L^2_\Pc (\Omega, \mathcal{F}_T, \mathbb{P};\mathbb{R}),
\end{align*}
Moreover $\hat{Y}^{x,a}  \in L^{2, loc}_\Pc(\Omega;{C} ([0,+\infty[;\mathbb{R}))$.
\newline Finally \eqref{ErgodicRan} is verified by $(\hat{Y}^{x,a},\hat{Z}^{x,a},\hat{K}^{x,a})$, $\mathbb{P}$-a.s. for all $t$ and $T$ with   $0\leq t\leq T$.
\end{Theorem}
\Dim
To simplify notation, we denote $(Y^{x,a,\beta},Z^{x,a,\beta},K^{x,a,\beta})$ simply by $(Y^\beta,Z^\beta,K^\beta)$.
\newline Set $\hat{Y}^\beta= \hat{v}^\beta (\X^{x,a}_t)= v^\beta (\X^{x,a}_t)- v ^\beta (0)$. By trivial computations, we find
\[
\hat{Y}^\beta_t \ = \ \hat{v}^\beta(\X^{x,a}_T)  -\beta \int_t^T \hat{Y}^\beta _s \, ds - \beta v^\beta(0)(T-t)+ \int_t^T \ell(\X^{x,a}_s, \I^{a}_s)\, ds - K^\beta_T + K ^\beta_t- \int_t^T Z^\beta_s  \, d W^1_s.
\]
From \eqref{CondEll} and \eqref{defvbar}, we have:
\begin{align*}
\lim_{h\to \infty}  \hat{Y}^{\beta_{k_{h}}}_t \ &= \ \hat {Y}^{x,a}_t, \qquad \forall\,t \in [0,T],\quad \mathbb{P}\hbox{-a.s.} \\
\lim_{h\to \infty}  \beta_{k_{h}} \hat{Y}^{\beta_{k_{h}}}_t  \ &= \ 0, \qquad\qquad\!\!\!\! \forall\,t \in [0,T], \quad \mathbb{P}\hbox{-a.s.}\\
\lim_{h\to \infty} \beta{_{k_{h}}}v^\beta_{k_{h}}(0) \ &= \ \lambda.
\end{align*}
Thanks to \eqref{CondZEll}, we know that there exists $\hat{Z}^{x,a}  \text { in } L^2_\Pc (\Omega\times [0,T];\Xi)$ such that
\begin{equation*}
Z^{\beta_{k_{h}}} \rightharpoonup \hat{Z}^{x,a}  \text { in } L^2_\Pc (\Omega\times [0,T];\Xi^*).
\end{equation*}
Thus all terms, apart from $K^\beta _T-K^\beta_t$, weakly converge. Choosing $t=0$,  we deduce that, for all $T\geq 0$, $K^\beta_T$, as well, converges weakly in $L^2(\mathcal{F}_T)$ to some $\hat{K}^{x,a}_T$. Therefore, we deduce 
that we can pass to the limit as $h\to \infty$ to get that $(\hat{Y},\hat{Z},\hat{K})$ is a solution to \eqref{ErgodicRan}.\finedim
\begin{Remark}
Notice that, thanks to  \eqref{unifLipsn}, we have
\[  |\hat{Y}_t| \leq L_\ell |X^{x,a}_t|, \qquad \forall\,t \in [0,T]. \]
\end{Remark}
We wish to prove that the above solution is maximal (in a suitable sense). 

The notion of maximality is the same as in \cite{CoFuhPham} but the proof is completely different since, once more, we can not use PDE interpretation. 

The following is, as a matter of fact, the main technical result of the paper.

To state it, fixed an horizon  $T>0$, letting $\lambda$ and $\hat{v}$ be the ones defined in \eqref{defvbar} and \eqref{deflambdabar}, we  consider the following BSDE on $[0,T]$
\begin{equation} \label{ErgodicRan-fino_a-T}
Y_t \ = \ \hat{v}(\X^{x,a}_T) + \int_t^T (\ell(\X^{x,a}_s,\I^a_s) -\lambda ) \,ds  - K_T + K_t - \int_t^T Z_s d W^1_s, \qquad  t\in[0,T].
\end{equation}
Notice that the above equation is indeed
  \eqref{EllitticaRandT}  with $\beta=0$, final condition $\phi$ replaced by $\hat{v}$ and generator $\ell$ replaced by $\ell-\lambda$,  see Definition \ref{SolEllRandT} for the definition of solution.
\begin{Theorem}\label{ergodico-massimale}
Assume {\bf(A.1)}--{\bf(A.8)}. The triple $(\hat{Y}^{x,a},\hat{Z}^{x,a},\hat{K}^{x,a})$  is the maximal solution of equation \eqref{ErgodicRan} in the sense that, for all fixed $T\geq 0$, if  $({Y},Z,K)$ defined on $[0,T]$ is a solution of \eqref{ErgodicRan-fino_a-T} then $\hat{Y}_t^{x,a} \geq Y_t$ for all $t \in [0,T]$. In particular, we have
\begin{equation}\label{corollario5.1}
\hat{v}(x) \ = \ \inf_{u \in  \mathcal{U}} \hat{J}^T (x,u), \qquad \forall\,x\in H,\,T>0,
\end{equation}
where
\[
\hat{J}^T (x,u) \ := \ \left[  \E \int_0^{T} (\ell (X^{x,u}_s,u_s) -\lambda) \,ds+  \hat{v}(X^{x,u}_T) \right].
\]

\end{Theorem}
\Dim
First of all, we know from Section \ref{sub-finhor}, Proposition \ref{propElliticaT}, that equation \eqref{ErgodicRan-fino_a-T}
admits a maximal solution  that can be rewritten as ${Y}_t=v^{0, T-t}(X^{x,a}_t) $ (for $v^{0, T-t}$ we refer to the notation in Proposition \ref{propElliticaT} with $\ell$ replaced by $\ell -\lambda$, we recall that $0$ indicates absence of discount, that is $\beta=0$). Moreover, by Remark \ref{SolEllRandT0} the following holds
\begin{equation}\label{sublinear} |v^{0,T-t}(X^{x,a}_t))| \leq c (1+ |X^{x,a}_t|), \qquad \forall\,t \in [0,T],\end{equation}
for some positive constant $c>0$, depending only  $T$ and on the constants introduced in Assumptions {\bf (A)}.

To prove that $v^{0,T-t}(X^{x,a}_t))\leq \hat{v}(X^{x,a}_t), \forall\,t \in [0,T]$, we come back to the penalized, finite horizon BSDE making a step ``backward''. As a matter of fact we introduce the classical (finite horizon) BSDE in $[0,T]$:
\begin{align*}
\check{Y}_t^{\beta,n} \ &= \ v^\beta(\X^{x,a}_T)
  -\beta \int_t^T \check{Y}_s^{\beta,n} \, ds + \int_t^T \ell(\X^{x,a}_s, \I^{a}_s)\,ds - n \int_t^T |\check{\Gamma}^{\beta,n}_s| \, ds \nonumber\\ 
&\quad \ -\int_t^T \check{Z}_s^{\beta,n} \, d W^1_s  -\int_t^T \check{\Gamma}^{\beta,n}_s \, d W^2_s,  \qquad 0 \leq t \leq T.
\end{align*}
See Proposition \ref{proprietan} for the definition on $v^{\beta}$. By Proposition \ref{propEllitica}, point \textit{4.}, we have that, for all $t \in [0,T]$,
\begin{equation*}
\check{Y}_t^{\beta,n} \searrow \check Y^\beta _t  \qquad \text{ as } n \to \infty,
\end{equation*}
where $\check{Y}_t^{\beta}$ is the  first component of the maximal solution to the following constrained BSDE on $[0,T]$:
$$
{Y}_t\ := \ v^\beta(\X^{x,a}_T)
  -\beta \int_t^T {Y}_s \, ds + \int_t^T \ell(\X^{x,a}_s, \I^{a}_s)\,ds -  K_T+K_t  -\int_t^T {Z}_s \, d W^1_s  \qquad 0 \leq t \leq T.
$$
Notice however that, by Proposition \ref{maximalellitica} the 
maximal solution of the above equation is indeed
$(Y^\beta,Z^\beta, K^\beta)$ (recall that we denote by $(Y^\beta,Z^\beta, K^\beta)$ the triple  $(Y^{x,a,\beta},Z^{x,a\beta}, K^{x,a\beta})$  founded on Proposition \ref{propEllitica}. So, in particular,
$
\check{Y}_t^{\beta,n} \searrow Y^\beta _t $ as  $n \to \infty$.
\newline Thus, setting $\hat{\check{Y}}_t^{\beta,n} = \check{Y}_t^{\beta,n} - v^\beta(0)$, we find,  for all $t \in [0,T]$: $\hat{\check{Y}}_t^{\beta,n} \searrow Y^\beta _t -v^\beta(0)=\hat{Y}^\beta _t$  as $ n \to \infty$.
\newline Moreover, for all $t\leq T$, it holds, $\mathbb{P}$-a.s.,
\begin{align*}
\hat{\check{Y}}_t^{\beta,n} \ &= \ \hat{v}^\beta(\X^{x,a}_T)
  -\beta \int_t^T \hat{\check{Y}}_s^{\beta,n} \, ds -\beta \int_t^T v^{\beta}(0) \, ds + \int_t^T \ell(\X^{x,a}_s, \I^{a}_s)\,ds - n \int_t^T |\hat{\check{\Gamma}}^{\beta,n}_s| \, ds \nonumber\\ 
&\quad \ -\int_t^T \hat{\check{Z}}_s^{\beta,n} \, d W^1_s  -\int_t^T \hat{\check{\Gamma}}^{\beta,n}_s \, d W^2_s.
\end{align*}
Then $( \hat{\check{Y}}^{\beta,n}, \hat{\check{Z}}^{\beta,n}, \hat{\check{\Gamma}}^{\beta,n})$ becomes the candidate triple to be compared with a generic solution $(Y,Z,K)$ of equation \eqref{ErgodicRan} such that $Y_T=\hat{v}(\X_T^{x,a})$ (recall that $\hat{Y}^{\beta_k}_t \rightarrow \hat{Y}^{\beta}_t$ as $k\rightarrow \infty$).
Notice that,  setting ${Y}^{\sharp}_t= e^{-\beta t} Y_t$, $Y_t^{\sharp,\beta,n} =  e^{-\beta t}\hat{\check{Y}}_t^{\beta,n}$, ${Z}^{\sharp}_t= e^{-\beta t} Z_t$, $Z_t^{\sharp,\beta,n} =  e^{-\beta t}\hat{\check{Z}}_t^{\beta,n}$, $\Gamma_t^{\sharp,\beta,n} =  e^{-\beta t}\hat{\check{\Gamma}}_t^{\beta,n}$, we find, $ \mathbb{P}$-a.s. for all $t\in [0,T]$:
\begin{align*} 
Y_t^{\sharp,\beta,n}-{Y}^{\sharp}_t \ = \  & e^{-\beta T}\hat{v}^\beta(\X^{x,a}_T) - e^{-\beta T}\hat{v}(\X^{x,a}_T) -   \int_t^T  e^{-\beta s}[ \beta v^\beta (0)-\lambda] \, ds+  \beta \int_t^T {Y}^{\sharp}_s\, ds \\ 
& - n \int_t^T  |\Gamma^{\sharp,\beta,n}_s| \, ds + \int_t^T e^{-\beta s} d\, K_s  + \int_t^T Z_s^{\sharp,\beta,n} \, d W^1_s  -\int_t^T \Gamma^{\sharp,\beta,n}_s \, d W^2_s. \notag
\end{align*}
By Girsanov's theorem there exists an equivalent probability $\tilde{\mathbb{P}}^{\beta,n}$ and a  $\tilde{\mathbb{P}}^{\beta,n}$-Wiener process $\tilde{W}^2$ such that ${W}^1$ is still a Wiener process under $\tilde{\mathbb{P}}^{\beta,n}$  and the above  equation
 can be rewritten as:
\begin{align*} \label{EllitticaRandDiff}
Y_t^{\sharp,\beta,n}-{Y}^{\sharp}_t \ = \  & e^{-\beta T}\hat{v}^\beta(\X^{x,a}_T) - e^{-\beta T}\hat{v}(\X^{x,a}_T) -   \int_t^T  e^{-\beta s}[ \beta v^\beta (0)-\lambda] \, ds+  \beta \int_t^T {Y}^{\sharp}_s\, ds \\ 
& + \int_t^T e^{-\beta s} d\, K_s  + \int_t^T Z_s^{\sharp,\beta,n} \, d W^1_s  -\int_t^T \Gamma^{\sharp,\beta,n}_s \, d W^2_s. \notag
\end{align*}
Hence, taking the expectation with respect to 
$\tilde{\mathbb{P}}^{\beta,n}$, we obtain
\begin{align*}
Y_0^{\sharp,\beta,n}-{Y}^{\sharp}_0  &\geq (\lambda - \beta v^\beta(0)) \frac{e^{-\beta T }-1}{\beta} + \beta \tilde{\E}^{\beta,n} \int_t^T {Y}^{\sharp}_s\, ds + e^{-\beta T}  \tilde {\E}^{\beta,n}  
[\hat{v}^\beta(\X^{x,a}_T) - \hat{v}(\X^{x,a}_T)] \\ & =   I^\beta_1+ I_2^{\beta,n}+I_3^{\beta,n}.
\end{align*} 
We have immediately that:
\[  \lim_{\beta \to \infty}    I^\beta_1 = 0.\]
Regarding the other two terms, we notice that the law of $\X^{x,a}$, under $\tilde{\mathbb{P}}^{\beta,n}$, is the law, under $\mathbb{P}$, of the solution
$\X^{x,a,\beta,n}$ of
\[
\begin{cases}
d\X_t\ = \ A\X_t dt + F(\X_t, \I^{\beta,n}_t)dt +\,GdW_t^1, \hspace{10mm} \X_0 \ = \ x, \\
d \I^{\beta,n}_t\ = \ R\gamma^{\beta,n}_t dt + RdW_t^2,  \hspace{32.3mm} \I_0 \ = \ a, \\
\end{cases}
\]
where 
\[
\gamma^{\beta,n}_s  = \left\{ \begin{array}{ll}- n  \frac{\Gamma^{\sharp,\beta,n}_s }{| \Gamma^{\sharp,\beta,n}_s |}, & \Gamma^{\sharp,\beta,n}_s \not= 0, \\
0, &   \text{ elsewhere. }\end{array} \right.
\]
Comparing with \eqref{System} we also notice that $ \X_t^{x,a,\beta,n}= \X_t^{x,a,\gamma^{\beta,n}}$. 
\newline By Remark \ref{R:J^R=J} and Proposition \ref{Prop-esunstate} we have that, for any  $p\geq 1$,
    \begin{equation}\label{stimaXp}\E  \sup_{t \in [0,T]} | \X_t^{x,a,\beta,n} |^p_{H}= \tilde{\E}^{\beta,n} \sup_{t \in [0,T]} | \X_t^{x,a} |^p_{H} \leq C, 
\end{equation}    
where the positive constant $C$ may depend on $x$, $T$, $p$, but neither on $n$ nor on $\beta$, $a$, $R$, $\gamma^{\beta, n}$ (since the dissipativity on $F$ is uniform with respect to the control variable).
Moreover, by \eqref{sublinear}, we get, for every $t\in [0,T]$,
\[  |{Y}^{\sharp}_t|  \leq C e^{-\beta t} (1+|\X_t^{x,a,\beta,n}|).\]
This together with \eqref{stimaXp} yields:
\[ \ \tilde{\E}^{\beta,n} \beta \int_0^T |{Y}^{\sharp}_t|\, dt \rightarrow 0, \qquad \hbox{ as $\beta \rightarrow 0$}, \]
that is
\[ \lim_{\beta \to 0} I^{\beta,n}_2 =0. \]
It remains to study the term $I_3^{\beta,n}$.
We notice that
\[ \tilde \E^{\beta,n} [\hat{v}^\beta(\X^{x,a}_T) - \hat{v}(\X^{x,a}_T)]  = \E [\hat{v}^\beta(\X^{x,a,\beta, n}_T) - \hat{v}(\X^{x,a,\beta, n}_T)].\]
Thanks to \eqref{unifLipsn}, we deduce that $(\hat{v}^\beta)_\beta$ is   equicontinuous and equibounded  on $H$.
Taking into account  \eqref{stima_compattezza}, we obtain
\[
\E || {\X^{x,a,\beta, n}_T}||_{D((\delta I-A)^{\rho})}\leq  C,
\]
where again $C$ depends on $T$ and on the quantities introduced in the Assumptions, but neither on $n$ nor on $\beta$.
Thus, for every $R>0$, set $B_R$ the centred ball of radius $R$ in  $D(\delta I -A)^{\rho}$ and denote by $B_R^c$ its complementary set. Then, we have
\[
|\E [\hat{v}^\beta(\X^{x,a,\beta, n}_T) - \hat{v}(\X^{x,a,\beta, n}_T)]| \leq  \ \sup_{y\in B_R}  |\hat{v}^\beta(y)-\hat{v}(y)|  + \E [|\hat{v}^\beta(\X^{x,a,\beta, n}_T) - \hat{v}(\X^{x,a,\beta, n}_T)| I_{B^c_R}(\X^{x,a,\beta, n}_T)].
\]
By Ascoli-Arzel\'a Theorem, $ \hat{v}^\beta$ converges uniformly on compact subsets of $H$ as $\beta$ tends to $0$, thus recalling \textbf{(A.3)}:
\[ \lim_{\beta \to 0}  \sup_{y\in B_R}  |\hat{v}^\beta(y)-\hat{v}(y)|=0. \]
On the other hand, \eqref{crescita_v} yields:
\begin{align*}
&\E [|\hat{v}^\beta(\X^{x,a,\beta, n}_T) - \hat{v}(\X^{x,a,\beta, n}_T)| I_{B^c_R}(X^{x,a,\beta, n}_T) ]\\ &\leq 
(\E [|\hat{v}^\beta(\X^{x,a,\beta, n}_T) - \hat{v}(\X^{x,a,\beta, n}_T)|^2])^{1/2}( \mathbb{P}(|| {\X^{x,a,\beta, n}_T}||_{D(\delta I-A)^{\rho})} >R))^{1/2} \\& 
\leq C( \E (|\X^{x,a,\beta, n}_T|^2))^{1/2}  (\E || {\X^{x,a,\beta, n}_T}||_{\mathcal{D}(\delta I-A)^\rho})^{1/2} R^{-1/2} \leq C R^{-1/2}, 
\end{align*}
for some positive constant $C$ independent of $n$ and $\beta$.
Therefore for every $\varepsilon >0$, we can find  $R$ large enough and then $\beta$ small enough such that
\[ |I^\beta_3| < \varepsilon. \]
This concludes the proof.
\finedim

\section{Long-time asymptotics of the finite horizon problem and the Ergodic control Problem}
\label{S:Finite}
Once the existence of a maximal solution $(\hat{Y}^{x,a},\hat{Z}^{x,a},\hat{K}^{x,a})$ and $\lambda$ to the constrained ergodic BSDE has been proved together with the Markovian representation of $\hat{Y}_t$ in terms of $\hat v(\X^{x,a})$, this final section is devoted to the application of such results to the study of the asymptotic expansion of the value function 
of finite horizon problems as $T\rightarrow +\infty$ as well as to the study of the value function of an ergodic optimal control problem.

We notice that the key point for the results below is that we have been able to identify through the constrained ergodic BSDE a function $\hat v$ (independent of $T$) and a constant $\lambda$ that, for all $T>0$, gives the value function of a control problem with horizon $T$ running cost $\ell -\lambda$ and final cost $\hat v$ itself.

\begin{Theorem}\label{T:LongTimeAsympt}
Assume {\bf(A.1)}--{\bf(A.5)} and {\bf(A.7)}. Fix $T>0$ and choose any function $\phi$ that satisfies {\bf (A.6)}. Let $v^T:=v^{0,T}$ be defined as in Proposition \ref{RegStatoU}  with $\beta=0$. For every $x\in  H$ it holds:
\[
|v^T(x) - \hat{v}(x) -\lambda T| \leq   C(1+|x|),
\]
for some positive constant $ C$. In particular
\[
\lim_{T\rightarrow+\infty} \frac{v^T(x)}{T} \ = \ \lambda,
\]
(in particular the unknown $\lambda$ in the ergodic BSDE \eqref{ErgodicRan} is uniquely determined).
\end{Theorem}

\Dim
By equality \eqref{corollario5.1} and Proposition \ref{RegStatoU} (with $\phi=\hat v$ and $\ell$ replaced by $\ell-\lambda$) we obtain

\begin{align*}
|v^T(x) - \hat{v}(x) - \lambda T  | \leq  \sup_{u \in  \mathcal{U}}  |J^T (x,u) - \hat{J}^T (x,u) -\lambda T| \leq \sup_{u \in  \mathcal{U}}  \big(\E |  \phi(X^{x,u}_T)| + \E |  \hat{v}(X^{x,u}_T)|\big). 
\end{align*}
By \ref{crescita_v} and {\bf (A.7)}  we deduce
\[ |  \hat{v}(X^{x,u}_T)| + | \phi(X^{x,u}_T)| \leq (1+|X^{x,u}_T|).  \]
and the claim follows by \eqref{stimadisstato-oriz-fin-prod}.
\finedim

\begin{Proposition}\label{P:LongTimeAsympt}
Assume  {\bf(A.1)}--{\bf(A.7)} and let $x\in H$.
\begin{itemize}
\item[\textup{(i)}] The real number $\lambda$ in \eqref{deflambdabar} satisfies the following inequality:
\[
\lambda \ \leq \ \inf_{u\in\Uc} \liminf_{T\rightarrow+\infty} \frac{1}{T}\E\bigg[\int_0^T \ell(X_t^{x,u},u_t)\,dt\bigg].
\]
\item[\textup{(ii)}] Assume in addition that there exists $\hat u\in\Uc$ such that the stochastic process
\[
\bigg(\hat v(X_t^{x,\hat u}) + \int_0^t \ell(X_s^{x,\hat u},\hat u_s)\,ds - \lambda\,t\bigg)_{t\geq0}
\]
is a local  martingale. Then
\[
\lambda \ = \ \lim_{T\rightarrow+\infty} \frac{1}{T}\E\bigg[\int_0^T \ell(X_t^{x,\hat{u}},\hat{u}_t)\,dt\bigg] \ \inf_{u\in\Uc} \liminf_{T\rightarrow+\infty} \frac{1}{T}\E\bigg[\int_0^T \ell(X_t^{x,u},u_t)\,dt\bigg].
\]
\end{itemize}
\end{Proposition}
\begin{Remark} We observe that the functional
$$\hat{J}(x,u):= \liminf_{T\rightarrow+\infty} \frac{1}{T}\E\bigg[\int_0^T \ell(X_t^{x,u},u_t)\,dt\bigg]$$
is a natural definition of \textit{ergodic cost}.
Thus the result of Proposition \ref{P:LongTimeAsympt} above can be rephrased saying that $lambda $ is the optimal cost for the ergodic cost $\hat{J}(x,\cdot)$ (independently on $x$) and $\hat{u}$ is the optimal  control)
\end{Remark}
\begin{Remark}\label{R:lambda}
Notice that a similar result to Proposition 5.1 was proved in the finite-dimensional context in \cite{CoFuhPham}, Remark 5.6, based on PDE, rather than probabilistic, techniques. More precisely, in \cite{CoFuhPham} the authors assume the existence of a (Lipschitz) feedback control which is, in a suitable sense, optimal for the ergodic PDE, while here in item (ii) above we assume the existence of a control $\hat u\in\Uc$ (not necessarily in feedback form) which is optimal according to the martingale principle of optimality.
\end{Remark}
\textbf{Proof.}
We split the proof into two steps.

\vspace{2mm}

\noindent\emph{Proof of \textup{(i)}.} Recall from Theorem \ref{T:LongTimeAsympt} that
\[
\lambda \ = \ \lim_{T\rightarrow+\infty} \frac{v^T(x)}{T}.
\]
By its definition $v^T(x)\leq\E[\int_0^T \ell(X_t^{x,u},u_t)\,dt+\phi(X_T^{x,u})]$ for any $u\in\Uc$, we obtain
\[
\lambda \ \leq \ \liminf_{T\rightarrow+\infty}\frac{1}{T}\E\bigg[\int_0^T \ell(X_t^{x,u},u_t)\,dt+\phi(X_T^{x,u})\bigg], \qquad \forall\,u\in\Uc.
\]
By Assumption {\bf (A.7)} and \eqref{stimadisstato-oriz-fin-prod}, we deduce that $\liminf_{T\rightarrow+\infty}\E[\phi(X_T^{x,u})]/T=0$. Consequently
\[
\lambda \ \leq \ \inf_{u\in\Uc}  \liminf_{T\rightarrow+\infty}\frac{1}{T}\E\bigg[\int_0^T \ell(X_t^{x,u},u_t)\,dt\bigg].
\]

\vspace{2mm}

\noindent\emph{Proof of \textup{(ii)}.} By point (i), it is enough to prove the inequality
\[
\lambda \ \geq \ \inf_{u\in\Uc} \liminf_{T\rightarrow+\infty} \frac{1}{T}\E\bigg[\int_0^T \ell(X_t^{x,u},u_t)\,dt\bigg].
\]
By the local martingale property of $(\hat v(X_t^{x,\hat u}) + \int_0^t \ell(X_s^{x,\hat u},\hat u_s)\,ds - \lambda\,t)_{t\geq0}$, it follows, through a straightforward localization argument, that, for every $T>0$,
\[
\hat v(x) \ = \ \E\bigg[\hat v(X_T^{x,\hat u}) + \int_0^T \ell(X_t^{x,\hat u},\hat u_t)\,dt - \lambda\,T\bigg],
\]
which can be rewritten as
\[
\lambda \ = \ \E\bigg[\frac{1}{T}\int_0^T \ell(X_t^{x,\hat u},\hat u_t)\,dt + \frac{\hat v(X_T^{x,\hat u}) - \hat v(x)}{T}\bigg].
\]
By \eqref{crescita_v} and \eqref{stimadisstato-oriz-infin-prod}, we have
\[
\lim_{T\rightarrow+\infty} \E\bigg[\frac{\hat v(X_T^{x,\hat u}) - \hat v(x)}{T}\bigg] \ = \ 0.
\]
In conclusion, we obtain
\[
\lambda \ = \ \lim_{T\rightarrow+\infty} \frac{1}{T}\E\bigg[\int_0^T \ell(X_t^{x,\hat u},\hat u_t)\,dt\bigg] \ \geq \ \inf_{u\in\Uc} \liminf_{T\rightarrow+\infty} \frac{1}{T}\E\bigg[\int_0^T \ell(X_t^{x,u},u_t)\,dt\bigg].
\]
\ep
\begin{Remark} The above results have been stated referring to the product space formulation of the problem, see Section \ref{SubS:InfHorProductSpace}. Thanks to the equality of the value functions, see Proposition \ref{UguaglianzaValueFunction1}, identical results can be proved with the same arguments in the original formulation of the control problem, see Section \ref{sub-originalcontrol}.
\end{Remark}
\begin{Example}{\rm
We consider an ergodic control problem for a stochastic heat equation as in \cite{DebHuTess}. The difference here is that we can handle the degenerate noise without the {\em structure condition}; in particular, in the equation below, the diffusion coefficient $\sigma$ need not to stay away from 0, on the other hand we need to ask strong dissipativity for the drift $f$. 
As a matter of fact we consider:
\[
\left\{  \begin{array}{ll}
d_t X^{u}(t, \xi)= [ \frac{\partial}{\partial \xi ^2}X^{u}(t, \xi) + f(\xi, X^{u}(t, \xi), u(t,\xi))  ] \, dt + \sigma (\xi) \dot{W}(t, \xi) \, dt, & t \geq 0, \ \xi \in [0,1], \\ \\
X^{u}(t,0)=X^{u}(t,1)= 0,  \\
X^{u}(0,\xi)= x_0(\xi),
\end{array}
\right.
\]
where $W$ is the state-space white noise  on $[0,+\infty) \times [0,1]$.  An admissible control $u$ is a predictable process $ u : \Omega \times [0,+\infty) \times [0,1] \to \R$, such that $ \int_0^1u ^2 (\omega, t,\xi)\,d \xi < \infty$, $\forall\,t \in[0,+\infty )$, $\P$-a.s..
The cost functional is
\[
J(x_0,u)= \liminf_{T \to +\infty} \, \frac{1}{T}\,  \E \int_0^T \int_0^1 \ell(t, X^{u}(t, \xi), u(t,\xi)) \, d \xi \, ds.
\]

The abstract formulation in $H=L^2(0,1)$ and $U=  L^2(0,1)$ follows as in \cite[section 5]{BismutElworthy}.
We notice that the realization of the second-order derivative with Dirichlet  boundary conditions fulfills  assumptions {\bf (A.1)} and {\bf (A.3)}, see \cite{Lunardi}.

Then, Theorem  \ref{T:LongTimeAsympt} and Proposition \ref{P:LongTimeAsympt} apply provided that one asks, for instance,
\begin{enumerate}
\item $f : [0,1] \times \R^2 \to \R$ is a continuous function such that
\begin{align*}
|f(\xi, x,u)| &\leq M_f ( 1+ |x| +|u|), \\
|f(\xi, x,u)-f(\xi, x',u) | &\leq  L_f |x-x'|,
\end{align*}
for suitable positive constants  $M_f$, $L_f$, for almost all $ \xi \in [0,1]$ and every $x, x', u \in \R$. Moreover we assume $f  (\xi, \cdot, u )\in  C^1(\R)$ for almost all $ \xi \in [0,1]$ and every $u \in \R$ to be such that, for some $\mu>0$,
\[ \frac{\partial}{\partial x}  f  (\xi, x, u ) \leq -\mu, \]
for almost all $ \xi \in [0,1]$ and every $x,u \in \R $.
\item  $ \sigma: [0,1] \to \R$ is a measurable and bounded function.
\item $\ell: [0,1] \times \R^2 \to \R$ is a continuous and bounded function such that
\[ |\ell(\xi, x,u)- \ell(\xi, x',u)| \leq L_\ell| x-x'|,
\]
 for a suitable positive constant $L_\ell$, for almost all $ \xi \in [0,1]$ and every $x, x', u \in \R$.
\end{enumerate}
}
\end{Example}
\begin{Example}{\rm
We can  also handle an SPDE, in a bounded open regular domain $D \subset \R^2$ driven by colored (in space) noise. We notice that in this case the structure condition would be a very artificial request. Namely we consider 
\[
\left\{  \begin{array}{ll}
d_t X^{u}(t, \xi)= [\mathcal{A}X^{u}(t, \xi) + f(\xi, X^{u}(t, \xi), u(t,\xi))  ] \, dt + \displaystyle{\frac{\partial{\mathcal{W}^G}}{\partial t}(t, \xi) \, dt}, & t \geq 0, \xi \in D, \\ \\
X^{u}(t,0)=X^{u}(t,1)= 0,  \\
X^{u}(0,\xi)= x_0(\xi),
\end{array}
\right.
\]
where

\begin{enumerate}
\item  $\mathcal{A}$ is a second-order operator in divergence form  $ \mathcal{A} =  \frac{\partial}{\partial \xi_h}\left( \sum_{h,k =1}^2 a_{h,k}(\xi) \frac{\partial}{\partial \xi_k}\right)$, the matrix $(a_{h,k})_{h,k}$ is non-negative and symmetric, with regular coefficients; moreover, the uniform elliptic condition
\[  \sup_{\xi \in \bar{D}}  \sum_{h,k =1}^2 a_{h,k}(\xi) \lambda_h \lambda_k \geq \mu |\lambda| ^2, \qquad \lambda \in \R^2. \]
is fulfilled for a positive constant $\mu$.

Under such an assumption the realization of $\mathcal{A}$ in $H=L^2(D)$ is a self-adjoint  strongly dissipative operator  $A$ that generates an analytic semigroup with dense domain in $H$, so that conditions {\bf (A.1)} and {\bf (A.3)} are fulfilled.

Moreover, it is  diagonalized by a complete orthonormal basis of functions $ e_{k}$ such that
\[
Ae_k = \lambda_k e_k,  \quad \text{ with }   \quad  \lambda_k  \sim  k, \qquad\qquad k =1,2,\dots
\]
\item $f : D \times \R^2 \to \R$ is a continuous function such that
\begin{align*}
|f(\xi, x,u)| &\leq M_f ( 1+ |x| +|u|), \\
|f(\xi, x,u)-f(\xi, x',u) | &\leq  L_f |x-x'|,
\end{align*}
for suitable positive constants  $M_f$, $L_f$, for almost all $ \xi \in D$ and every $x, x', u \in \R$. Moreover we assume $f  (\xi, \cdot, u )\in  C^1(\R)$ for almost all $ \xi \in D$ and every $u \in \R$; in addition, for some $\nu>0$,
\[ \frac{\partial}{\partial x}  f  (\xi, x, u ) \leq -\nu, \]
for almost all $ \xi \in D$ and every $x,u \in \R $.
\item $\ell: [0,1] \times \R^2 \to \R$ is a continuous and bounded function such that
\[ |\ell(\xi, x,u)- \ell(\xi, x',u)| \leq L_\ell| x-x'|,
\]
 for a suitable positive constant $L_\ell$, for almost all $ \xi \in [0,1]$ and every $x, x', u \in \R$.
\item An admissible control $u$ is a predictable process $ u : \Omega \times [0,+\infty) \times D \to \R$, such that \[ \int_D u ^2 (\omega, t,\xi)\,d \xi < \infty, \quad \forall\,t \in[0,+\infty ), \P\text{-a.s.}.\]
\item $\frac{\partial{\mathcal{W}^G}}{\partial t}(t, \xi)$ stands for a gaussian noise that is assumed to be white in time and colored in space. More precisely, the infinite dimensional reformulation of the equation  is driven  by an $H$-valued Wiener process $W^G$ defined using the  sum:
\[   W^G (t, \xi) = \sum_{k=1}^\infty G e_k(\xi) \beta_{k}(t),   \quad t \geq 0, \quad \xi \in D, \]
where $ \{ \beta_k\}$ is a sequence of mutually independent standard Brownian motions and  $G: H \to H$ is given by $G:= (-A)^{-\eta}$ for a fixed $\eta  >\frac{1}{4}$, in order to ensure that condition {\bf(A.2)} holds true. The operator $G$ is clearly not invertible.

Under the above assumptions, Theorem  \ref{T:LongTimeAsympt} and Proposition \ref{P:LongTimeAsympt} apply.
\end{enumerate}}
\end{Example}

\bibliographystyle{plain}
\bibliography{ERSubmitted}

\end{document}